\magnification\magstep1
\baselineskip = 18pt
\overfullrule = 0pt

\def\n{\noindent}
 \def\qed{{\hfill{\vrule height7pt width7pt
depth0pt}\par\bigskip}} 

\def\vp{\varepsilon}
\def\pf{\medskip\n {\bf Proof.}~~}

\def\ms{\medskip}

\def\ovl{\overline}

\def\FF{{\rm I}\!{\rm F}}
\def\CC{ \;Ê{}^{ {}_\vert }\!\!\!{\rm C}}

\def\RR{{\rm I}\!{\rm R}}

\def\tilde{\widetilde}
\def\hat{\widehat}
\def\ie{{\it i.e.\/}\ }
\def\cf{{\it cf.\/}\ } 

\input mssymb
\centerline{\bf Grothendieck's Theorem}
\centerline{\bf for Operator Spaces}\bigskip
\centerline{by Gilles Pisier\footnote*
{Partially supported by   NSF      and
   Texas Advanced Research
 Program 010366-163}}
\centerline{Texas A\&M University
 and Universit\'e Paris 6}
\centerline{and }\centerline{Dimitri
Shlyakhtenko\footnote{**}{Partially supported by
  NSF and  Sloan fellowship}}
 \centerline{University of California at
Los Angeles}\bigskip
\bigskip
\centerline{\bf Abstract} We prove several versions of
Grothendieck's Theorem for completely bounded
linear maps $T\colon \ E \to F^*$, when
$E$ and $F$ are operator spaces. We prove that
if $E,F$ are $C^*$-algebras,   of which at least one is exact,
then every completely bounded $T\colon \ E \to F^*$
can be factorized through the direct sum
of the row and column Hilbert operator spaces.
Equivalently $T$ can be decomposed
as $T=T_r+T_c$ where $T_r$ (resp. $T_c$) factors
completely boundedly through 
a row (resp. column) Hilbert operator space.
This settles positively (at least partially)
 some earlier
conjectures of Effros-Ruan and Blecher on
the factorization of completely bounded  bilinear
forms on $C^*$-algebras.
Moreover, our result holds more generally
for any pair $E,F$ of ``exact" operator spaces.
This  yields a characterization of the
completely bounded maps from a $C^*$-algebra
(or from an   exact operator space)
to the operator Hilbert space OH.
As a corollary we prove that,
up to a complete isomorphism, the row and column Hilbert operator spaces
and their direct sums are the only 
operator spaces $E$ such that both $E$ and its dual $E^*$ are exact.
We also characterize the  Schur multipliers which are completely bounded
from the space of compact operators to the trace class.

\medskip

\vfill\eject

\n {\bf Introduction.}
In 1956, Grothendieck published the fascinating paper [G]
now often referred to   as ``the R\'esum\'e".
The central result there was described by  Grothendieck
as ``the fundamental theorem of the metric theory
of tensor products". This result, now known as Grothendieck's theorem
  (GT in short) -or sometimes ``Grothendieck's inequality"- has played a 
major role in the developments of Banach space theory in the last
three decades; moreover, its non-commutative version has also found
important applications to several specific questions 
in $C^*$-algebra theory (see [P5]). It is   natural to wonder whether 
  this result still holds, at least in some form, for the 
recently introduced and currently very active non-commutative
analogue of Banach spaces, namely  ``operator spaces", in the sense
of [BP, ER2] (see also  [ER1, P1]), and this is precisely
  the goal of the present paper.
Let $A,B$ be $C^*$-algebras.
 While the previous versions
are all concerned with {\it bounded} bilinear forms
on $A\times B$,
or equivalently bounded  linear maps
$T\colon\ A\to B^*$, 
 and their possible {\it bounded} factorizations, we  will study
{\it completely bounded} bilinear forms
or equivalently completely bounded  linear maps
$T\colon\ A\to B^*$,
and their possible {\it completely bounded} factorizations.
For instance,   the classical GT and its later extensions say that any bounded
$T\colon\ A\to B^*$ factors boundedly through a Hilbert space.
One of our main results says that, under a mild restriction on either $A$ or $T$,
any completely bounded
$T\colon\ A\to B^*$ factors completely boundedly through the direct sum
of two very simple ``building blocks": the row and column Hilbert operator
spaces.
 We can thus claim that while GT
entirely elucidates the bounded case, we
analogously elucidate the completely  bounded one.

We will now describe more precisely 
the connections of our work with the existing literature and
the  conjectures which motivated 
it.

\n The non-commutative version of Grothendieck's
theorem says that, if $A,B$ are
$C^*$-algebras any bounded bilinear form $u\colon
\ A\times B\to
\CC$ satisfies the following inequality:\ for any finite sequence $(a_i,b_i)$
in $A\times B$ we have
 $$ {\left|\sum u(a_i,b_i)\right| \le K\|u\| \left\{\left\|\sum
a^*_ia_i\right\| + \left\| \sum a_ia^*_i\right\|\right\}^{1/2}
\left\{\left\|\sum b^*_ib_i\right\| + \left\|\sum
b_ib^*_i\right\|\right\}^{1/2}}\leqno(0.1)  $$
where $K$ is a numerical constant (independent of $A,B$ and $u$).

\n Moreover, by a rather simple application of the Hahn-Banach theorem, (0.1)
implies that there are states $f_1,f_2$ on $A$ and $g_1,g_2$ on $B$ such that
$$\forall(a,b)\in A\times B\qquad |u(a,b)| 
\le K\|u\| \{f_1(aa^*) +f_2(a^*a)
\}^{1/2} \{g_2(bb^*) +g_1(b^*b) \}^{1/2}.\leqno (0.2)$$
This was proved for commutative $C^*$-algebras by Grothendieck. The
non-commutative case was obtained in [P4] with an approximability assumption
and in [H1] in full generality. It is easy to deduce from (0.2)
that, assuming $\|u\|\le 1$ for simplicity, there is a decomposition
$u=u_1+u_2+u_3+u_4$ with
$$|u_1(a,b)|\le K  \{f_1(aa^*)  
\}^{1/2} \{ g_1(b^*b) \}^{1/2}\quad
 |u_2(a,b)|\le K  \{f_2(a^*a)  
\}^{1/2} \{ g_2(bb^*) \}^{1/2}$$
$$|u_3(a,b)|\le K  \{f_1(aa^*)  
\}^{1/2} \{ g_2(bb^*) \}^{1/2}\quad
 |u_4(a,b)|\le K  \{f_2(a^*a)  
\}^{1/2} \{ g_1(b^*b) \}^{1/2}.$$

In view of the recent development of operator space theory, it is natural to
look for a version of this theorem for {\sl  jointly completely bounded\/}
(j.c.b. in short) bilinear forms.
But here the terminology poses a   problem,  
there are two different notions in the bilinear case: the joint 
complete boundedness ([ER1-2  ,BP]) and 
  the complete boundedness in Christensen and Sinclair's sense ([CS1-2,
CES]). The second notion came first, was then called
simply ``complete boundedness" and has proved extremely fruitful. 
 To (reluctantly) conform with the already established tradition,
we will call the  first ones ``jointly c.b." although 
calling them c.b. would definitely be more natural  
from the viewpoint of operator space theory. Indeed, given a $C^*$-algebra or
more generally an operator space
$F\subset B(H)$, the latter theory provides us with a natural realization of the
dual $F^*$ as an operator space, so that we have a specific isometric
embedding $F^*\subset B({\cal H})$. Thus, given two 
$C^*$-algebras $A,B$ and operator spaces
$E\subset A$ and $F\subset B$, we  say that a
bilinear form
$u\colon
\ E\times F\to \CC$ is jointly completely bounded (in short j.c.b.) if the
associated linear map $\tilde u\colon \ E\to
F^*\subset B(\cal H)$ is c.b. Moreover we let
$$\|u\|_{jcb} = \|\tilde u\|_{cb(E,F^*)}.\leqno (0.3)$$
Using tensor products, this definition can be extended to the case
of bilinear forms $u$ with values in $B({\cal H })$, the preceding case
then corresponds to $\dim  (  {\cal H }) =1$. See \S 1 for   details.

Now  a bilinear map $u\colon\
E\times F\to \CC$ is called c.b.\ (in Christensen and Sinclair's sense) if the
bilinear forms
$$u_n\colon \ M_n(E)\times M_n(F)\to M_n$$ $$
((a_{ij}), (b_{ij})) \to \left(\sum\nolimits_k u(a_{ik},b_{kj})\right)_{ij}$$
are uniformly bounded, and by definition
$$\|u\|_{cb} = \sup_{n\ge 1} \|u_n\|.\leqno (0.4)$$
Here again, this definition can be extended to the case
of bilinear forms $u$ with values in $B({\cal H })$. 
See \S 1 for   more on this.

This notion is much better understood than the preceding one. In particular,
it is easy to show that a bilinear form $u\colon \ E\times F\to \CC$ is c.b.\
with $\|u\|_{cb} \le 1$ iff for all finite sequence $(a_i,b_i)$ in $E\times F$
we have
$$\left|\sum u(a_i,b_i)\right|\le \left\|\sum a_ia^*_i\right\|^{1/2} \left\|
\sum b^*_ib_i\right\|^{1/2}.
\leqno (0.5)$$
Moreover, assuming $E\subset A$ and $F\subset B$, then (0.5) implies that
there are states $f_1,g_1$ on $A$ and $B$ respectively such that
$$\forall(a,b)
\in E\times F\qquad\quad |u(a,b)| \le (f_1(aa^*) g_1(b^*b))^{1/2}.\leqno (0.5)'
$$ 
(Note that conversely (0.5)$' \Rightarrow$ (0.5) so (0.5) and (0.5)$'$ are
essentially equivalent.) Moreover there is
an extension of $u$ to $A\times B$ with the same c.b.\ norm.

When $u\colon \ E\times F\to \CC$ is c.b.\ the associated linear map $\tilde
u\colon \ E\to F^*$ admits a factorization of the form $E {\buildrel v\over
\longrightarrow} H_r {\buildrel w\over \longrightarrow} F^*$ through a row
Hilbert space with $\tilde u = wv$ and $\|v\|_{cb} \|w\|_{cb} = \|u\|_{cb}$
([ER3]). In particular, we have
$$\|\tilde u\|_{cb} \le \|u\|_{cb}$$
or equivalently
$$\|u\|_{jcb}\le \|u\|_{cb}.\leqno (0.6)$$
It is worthwhile to observe that the Christensen-Sinclair notion is not
symmetric:\hfill\break if $u\colon \ E\times F\to
\CC$ is c.b.\ the transposed bilinear form
${}^tu\colon \ F\times E\to \CC$ is not
necessarily c.b., while for j.c.b.\ forms this is
true and we do have
$$\|{}^tu\|_{jcb} = \|u\|_{jcb}.$$
(This is related to the basic fact from operator space theory that $\|\tilde
u\|_{cb} = \|\tilde u{}^*\|_{cb}$.) Thus (0.6) implies for any bilinear from
$v\colon \ E\times F\to \CC$
$$\|v\|_{jcb}\le \|{}^tv\|_{cb}\leqno (0.7)$$
and more generally (0.6) and (0.7) together yield
$$\|u+v\|_{jcb} \le \|u\|_{cb} + \|{}^tv\|_{cb}.\leqno (0.8)$$
We note that (see (0.5) above) $\|{}^tv\|_{cb}\le 1$
iff for any finite sequence $(a_i,b_i)$ in
$E\times F$ we have
$$\left|\sum v(a_i,b_i)\right| \le \left\|\sum a^*_ia_i\right\|^{1/2}
\left\|\sum b_ib^*_i\right\|^{1/2},\leqno (0.9)$$
or equivalently iff there are states $f_2,g_2$ on $A$ and $B$ respectively
such that
$$\forall(a,b)\in
E\times F\qquad |v(a,b)|\le (f_2(a^*a) g_2(bb^*)))^{1/2}.\leqno (0.9)'$$

Thus if $U = u+v$ with $\|u\|_{cb}\le 1$ and $\|{}^tv\|_{cb}\le 1$ we find
states $f_1,g_1,f_2,g_2$ such that
$$\forall(a,b) \in E\times F\qquad |U(a,b)|\le (f_1(aa^*) g_1(b^*b))^{1/2} +
(f_2(a^*a) g_2(bb^*))^{1/2}\leqno (0.10)$$
or equivalently for all finite sequences $(a_i,b_i)\in E\times F$ and all
$\lambda_i>0$
$$\left|\sum U(a_i,b_i)\right|\le \left\|\sum a_ia^*_i\right\|^{1/2}
\left\|\sum b^*_ib_i\right\|^{1/2} + \left\|\sum \lambda_i
a^*_ia_i\right\|^{1/2} \left\|\sum \lambda^{-1}_i b_ib^*_i\right\|^{1/2}.
\leqno (0.10)'$$
Conversely it can be shown that (0.10) implies the existence of decomposition
$U = u+v$ with $\max\{\|u\|_{cb}, \|{}^tv\|_{cb}\} \le 1$. (The proof of this
converse is less obvious than may seem at first glance. We give the details in
\S 2 below. We show there that (0.10) implies a decomposition $U = u+v$
together with states $f'_1,f'_2,g'_1,g'_2$ such that
$$|u(a,b)| \le (f'_1(aa^*) g'_1(b^*b))^{1/2} \hbox{ and } |v(a,b)| \le
(f'_2(a^*a) g'_2(bb^*))^{1/2}$$
but we apparently cannot do this with the original states $f_1,f_2,g_1,g_2$.)
\medskip
We now return to the above decomposition $u=u_1+u_2+u_3+u_4$. There 
$u_1$ and
$u_2$  are jointly  c.b., but $u_3$ and $u_4$ in
general are not. More precisely, $u_1+u_2$ is
jointly c.b. on $A\times B$, while $u_3+u_4$
is jointly c.b. on $A\times B^{op}$, where $B^{op}$ denotes the opposite
$C^*$-algebra (anti-isomorphic to $B$).
Probably led by similar observations, Effros and Ruan 
formulated in [ER2] (with $K=1$)  the following
\medskip

\proclaim Conjecture 0.1. Let $A,B$ be $C^*$-algebras and let $u\colon \
A\times B\to\CC$ be a j.c.b.\ bilinear form. Then there exist states
$f_1,f_2,g_1,g_2$ on $A,B$ respectively such that
$$|u(a,b)|\le K\|u\|_{jcb} (f_1(aa^*) g_1(b^*b))^{1/2} + (f_2(a^*a)
g_2(bb^*))^{1/2}\leqno \forall (a,b)\in A\times B$$
where $K$ is a numerical constant. 

Independently, Blecher [B1] was led to a similar
conjecture:

\proclaim Conjecture 0.2. Let $A=B({\cal H })$ and $B=B({\cal K })$ where
${\cal H,K}$ are Hilbert spaces.
There is a constant $K$ such that
for any $w$ in the algebraic tensor product $A\otimes B$
we have
$$\|w\| _{\wedge} \le K\max\{ \|w\|_h, \|{}^t w\|_h\}$$
where $\|w\| _{\wedge}$ is the operator space version of
the projective norm and where $\|w\|_h$ (resp. $\|{}^t w\|_h$) is the norm
in the Haagerup tensor product $A\otimes_h B$ (resp.  $B\otimes_h A$).
%\leqno(0.B)

In  the unpublished problem book of a 1993 conference, he also
asked the same questions when $A,B$ are $C^*$-algebras.
Moreover although he did not make it precise, Blecher implicitly conjectured
that the corresponding matricial norms were uniformly equivalent,  
so that the associated tensor products $A\otimes^{\wedge} B$
and $[A\otimes_h B] \cap {}^t [B\otimes_h A ]$ should be completely isomorphic.
At the time of this writing, we are unable to prove this (even when $A$
and $B$ are both commutative!).

  By [BP, ER2] we know that if $A,B$ are {\it arbitrary operator spaces}, we have
$$ \|w\| _{\wedge}= \sup\{ |\langle w, U  \rangle | \}\leqno(0.11)$$
where the supremum runs over all bilinear forms
$U\colon\ A\times B\to \CC $ with $\|U\|_{jcb}\le1$, and also
$$ \|w\| _{A\otimes_h B}= \sup\{ |\langle w, u  \rangle | \}\leqno(0.12)$$
where the supremum runs over all bilinear forms
$u\colon\ A\times B\to \CC $ with $\|u\|_{cb}\le1$. Thus Conjecture 0.2
  is equivalent by duality to the following

\proclaim Conjecture 0.2'. There is a constant $K$ such that  any
j.c.b.\ bilinear form $U\colon \ A\times B\to \CC$ can be decomposed as a
sum $U = u+v$ with 
$\|u\|_{cb} + \|{}^tv\|_{cb} \le K\|U\|_{jcb}$.

In other words Blecher conjectured that the estimate (0.8)
can be reversed on $A\times B$,
thus establishing a very nice (and simple)
relationship between the two notions of c.b.\ for
bilinear forms.

Our main objective in this paper is to   prove versions of these conjectures.
More precisely, we will prove Conjectures~0.1 and 0.2'
 assuming either $A$ or $B$ is exact,  or for
arbitrary $C^*$-algebras but under a suitable
approximability assumption on the bilinear form
$u$. Actually, our results are valid for j.c.b.\
bilinear forms
$U\colon
\ E\times F \to \CC$ defined on {\sl exact\/} operator spaces. In particular,
we will prove that, assuming $E\subset A$ and
$F\subset B$, any such map extends to a j.c.b.\
bilinear form
$\hat U$ defined on $A\times B$.
Unfortunately, since $B(H)$ is not exact, we cannot prove
Conjecture 0.2 as stated above. Also, the extension of
Conjecture 0.2' to the  case of a $B(H)$-valued bilinear form
 remains open at the time of this writing.

 Our results  
constitute a sequel to the papers [JP] and [P3].

\n We say that an operator space $E$ is exact  if there is a constant $C$
such that for any finite dimensional subspace $G\subset E$ there is an integer
$N$, a subspace $\widetilde G\subset M_N$ and an isomorphism $u\colon \ G\to
\widetilde G$ such that $\|u\|_{cb}\|u^{-1}\|_{cb}\le C$. We will denote by
$ex(E)$ the smallest $C$ for which this holds. This is the operator space analog
of a notion introduced and extensively studied by Kirchberg  for $C^*$-algebras  ([Ki3]).

In [JP], a characterization and an extension theorem were obtained for
``tracially'' bounded bilinear forms $u\colon \ E\times F\to \CC$ when $E,F$
are exact operator spaces. ``Tracial boundedness'' is a notion intermediate
between boundedness and complete boundedness. 

The precise statements of our main results  are as follows.

\proclaim Theorem 0.3. Let $E\subset A$, $F\subset B$ be exact operator spaces
sitting in $C^*$-algebras $A,B$ with exactness constants respectively $ex(E)$
and $ex(F)$. Let $C= ex(E) ex(F)$.
Then any j.c.b.\ bilinear form $U\colon \ E\times F\to \CC$ satisfies
the following inequality:\ for any finite sequence $(a_i,b_i)$ in $E\times F$
and for any $\lambda_i>0$ we have
$$ 
\left|\sum U(a_i,b_i)\right| \le \leqno (0.13) $$
 $$
 C\|u\|_{jcb} \left[\left\|
\sum \lambda_ia^*_ia_i\right\|^{1/2}\right. 
   \left. + \left\|\sum \lambda^{-1}_ia_ia^*_i
\right\|^{1/2}\right]
 \left[\left\|\sum
\lambda_ib^*_ib_i\right\|^{1/2} +
\left\|\sum \lambda^{-1}_ib_ib^*_i\right\|^{1/2}\right]
$$ and  consequently:
$$\left|\sum U(a_i,b_i)\right| 
\le
2 C\|u\|_{jcb} \left[\left\|
\sum a_ia^*_i\right\|^{1/2} \left\|\sum
b^*_ib_i\right\|^{1/2}\right. 
 \quad \left. +
\left\|\sum a^*_ia_i\right\|^{1/2} \left\|\sum
b_ib^*_i\right\|^{1/2}\right].\leqno(0.14)$$

\proclaim Theorem 0.4. Let $K=2^{3/2} C$.
Then for any j.c.b.\ bilinear form $U\colon \ E\times F\to \CC$
with $\|U\|_{jcb}\le 1$,  there are
states $f_1,f_2$ and $g_1,g_2$ on $A,B$ respectively such that
$$\forall(a,b)\in E\times F\qquad |U(a,b)| \le K   [(f_1(aa^*)
g_1(b^*b))^{1/2} + (f_2(a^*a) g_2(bb^*))^{1/2}].\leqno (0.15)$$
Moreover, for any finite sequence $(a_i,b_i)$
 in $E\times F$ and for any $\lambda_i>0$ we have
$$ {|\sum U(a_i,b_i)|   \le K \left[
\left\|\sum a_ia^*_i\right\|^{1/2} \left\|\sum b^*_ib_i\right\|^{1/2} 
   + \left\|\sum \lambda_ia^*_ia_i\right\|^{1/2} \left\|\sum
\lambda^{-1}_i b_ib^*_i\right\|^{1/2}\right].}\leqno (0.16)$$
Finally, any form $U$ satisfying (0.15) for some $K$ can be decomposed as
$U=u+v$ where $u,v$ are bilinear forms
satisfying
$$\max\{ \|u\|_h ,\|{}^t v\|_h\} \le K.$$

\proclaim Theorem 0.5. Let $A,B$ be $C^*$-algebras. Let
 $U\colon \
A\times B\to \CC$ be a bilinear form  with $\|U||_{jcb}\le 1$.
Either one of the two following assumptions ensures that 
$U$ satisfies the conclusions of Theorem~0.4 with $K =
2^{3/2}$.\hfil\break
\n (i) At least one of the algebras $A,B$ is exact.
\hfil\break
\n (ii) The form
   $U$ is the pointwise limit of a net of finite rank forms $U_\alpha$
such that $\|U_\alpha||_{jcb}\le 1$.

\proclaim Corollary 0.6. In the same situation as Theorem 0.4,
 any j.c.b.\ bilinear form
$U\colon \ E\times F\to \CC$ admits an extension $\hat U\colon \ A\times B\to
\CC$ with $\|\hat U\|_{jcb}  \le K\|U\|_{jcb}$.

Our result admits a nice reformulation as a
factorization involving the row and column
operator space structures (resp.\ $H_r$ and $H_c$)
on a Hilbert space
$H$:

\proclaim Corollary 0.7. In the   situation of Theorem 0.4, any
c.b.\ linear map $T\colon \ E\to F^*$ can
be factorized completely boundedly through an
operator space of the form $H_r\oplus H_c$ for
some Hilbert space $H$. More precisely there are
c.b.\ maps $v\colon \ E\to H_r\oplus H_c$ and
$w\colon \ H_r\oplus H_c\to F^*$ such that
$T = wv$ and $\|w\|_{cb} \|v\|_{cb} \le
2^{3/2} C\|T\|_{cb}$.
On the other hand, in the situation of
Theorem 0.5, any c.b.\ linear map $T\colon \ A\to B^*$
which can be approximated by finite rank   completely contractive maps,
 in the point weak-$*$ topology, 
 can
be factorized in the same way, $T = wv$ with $\|w\|_{cb} \|v\|_{cb} \le
2^{3/2}$; moreover,
 the approximability assumption on   $T$ can be dispensed with
if either $A$ or $B$ is exact.

In the ``predual'' situation, our results also yield:

\proclaim Corollary 0.8. Let $E,F$ be exact
 operator spaces and let $C = ex(E)
ex(F)$ as before. Let $A,B$ be arbitrary $C^*$-algebras. Then the following
isomorphisms hold:\hfil\break
\n (i) $E\otimes^{\wedge}  F\simeq [E\otimes _h
F]
\cap [{}^t(F \otimes_h E)]$. Here
$E\otimes^{\wedge}  F$ denotes the operator space
version of the projective tensor
product.\hfil\break
\n (ii)\ $A\otimes^{\wedge} B\simeq [A\otimes _h B]
\cap [{}^t(B \otimes_h A)]$.\hfil\break
\n (iii)\ $A^* \otimes_{\rm min} B^*\simeq A^* \otimes_h B^* + {}^t(B^*
\otimes_h A^*)$. \hfil\break The equivalence constant is $\le
4\sqrt{2}\ C$ in (i) and $\le 4\sqrt{2} $ in (ii) and (iii).

\n {\bf \S 1. Background}

\n We refer to [Pa1] (resp.\ [CS1-2]) for background on c.b. linear (resp.\
multilinear) maps and to [ER1] and [P1] for background on operator spaces in
general. We just recall that the minimal tensor product $E \otimes_{\rm min}
F$ of two operator spaces $E\subset B(H_1)$, $F\subset B(H_2)$ is defined so
that we have a completely isometric embedding
$$E \otimes_{\rm min} F\subset B(H_1 \otimes_2 H_2).$$
We now give a precise
 definition for ``joint'' complete boundedness.

\proclaim Definition 1.1. Let $E,F,G$ be operator spaces.  A bilinear form
$u\colon \ E\times F\to G$ is j.c.b.\ if for any
$C^*$-algebras $B_1,B_2$ $u$ can be ``extended''
to a bounded bilinear form $(u)_{B_1,B_2}\colon \
E\otimes_{\rm min} B_1\times F\otimes_{\rm min}
B_2 \to G\otimes_{\rm min} B_1
\otimes_{\rm min} B_2$ taking $(e\otimes b_1, f\otimes b_2)$ to $u(e,f)
\otimes b_1 \otimes b_2$. Moreover we have
$$\|u\|_{jcb} = \sup\|(u)_{B_1,B_2}\|\leqno (1.1)$$
where the supremum runs over all pairs $B_1,B_2$ of $C^*$-algebras.

Equivalently, we have $\|u\|_{jcb}\le 1$ iff for all $C^*$-algebras $B_1,B_2$
and for all finite sums
$$\sum a_i \otimes x_i \in E\otimes B_1,\quad \sum b_j\otimes y_j \in
F\otimes B_2$$
we have
$$\left\|\sum u(a_i,b_j) \otimes x_i\otimes y_j\right\|_{\rm min} \le
\left\|\sum a_i\otimes x_i\right\|_{\rm min} \left\|\sum b_j\otimes
y_j\right\|_{\rm min}.\leqno (1.1)'$$
In addition, it suffices to consider matricial $C^*$-algebras for $B_1$ and
$B_2$, more precisely we have actually
$$\|u\|_{jcb} = \sup_{n\ge 1} \|(u)_{M_n,M_n}\|.\leqno (1.2)$$
 It is
easy to see that the usual definition of $\|u\|_{jcb}$ (which is $\|u\|_{jcb}
= \|\tilde u\|_{cb(E,F^*)}$) is equivalent to (1.2).
The equality of (1.1) and (1.2) is a routine
verification left to the reader.

Following the same pattern we have the following equivalent definition for
complete boundedness:

\proclaim Proposition 1.2. A bilinear
 form $u\colon \ E\times F\to G$ is c.b.\
iff for any $C^*$-algebra $A$, $u$ can be ``extended'' to a bounded bilinear
form $u_A\colon \ E\otimes_{\rm min} A\times F\otimes_{\rm min} A\to
G\otimes_{\rm min} A$ taking $(e\otimes a_1, f\otimes a_2)$ to $u(e,f) \otimes
a_1a_2$. Moreover we have
$$\|u\|_{cb} = \sup\|u_A\|$$
where the supremum runs over all possible $C^*$-algebras $A$ and actually
$$\|u\|_{cb} = \sup_{n\ge 1} \|u_{M_n}\|.$$

The equivalences (0.1) $\Leftrightarrow (0.1)'$,
(0.5)~$\Leftrightarrow (0.5)'$, (0.9) $\Leftrightarrow~(0.9)'$ all follow
from the same Hahn-Banach type argument as follows.

\proclaim Proposition 1.3. Let $A,B$ be $C^*$-algebras. Assume $E\subset A$,
$F\subset B$. Let $u\colon \ E\times F\to \CC$ be a bilinear form. Let
$\alpha_1,\alpha_2$ and $\theta_1,\theta_2$ be fixed non-negative numbers.
Then the following assertions are equivalent
\item{(i)} There are states $(f_1,f_2), (g_1,g_2)$ on $A$ and $B$ respectively
such that for all $(a,b)$ in $E\times F$ we have 
$$|u(a,b)| \le [\alpha_1f_1(a^*a) + \alpha_2f_2(aa^*)]^{1/2}
[\theta_2g_2(b^*b) + \theta_1g_1(bb^*)]^{1/2}.$$
\item{(ii)} For all finite sequences $(a_i,b_i)$ in $E\times F$ we have
$$\left|\sum u(a_i,b_i)\right|\le \left[\alpha_1\left\|\sum a^*_ia_i\right\| +
\alpha_2 \left\|\sum a_ia^*_i\right\|\right]^{1/2} \left[\theta_2\left\| \sum
b^*_ib_i\right\| + \theta_1\left\|\sum b_ib^*_i\right\|\right]^{1/2}.$$

\pf
(sketch) (i) $\Rightarrow$ (ii) is obvious by Cauchy-Schwarz. Conversely
assume (ii). We will use the classical arithmetic/geometric mean inequality as
follows:
$$\forall~x,y\ge 0\qquad (xy)^{1/2} \le 2^{-1}(x+y) \hbox{ and }
(xy)^{1/2} = \inf_{\lambda>0} 2^{-1}(\lambda x + \lambda^{-1}y).\leqno
(1.3)$$ 
Thus (ii) and (1.3) imply: $$2\left|\sum u(a_i,b_i)\right|\le \sup\left[
\alpha_1\sum f_1(a^*_ia_i) + \alpha_2 \sum f_2(a_ia^*_i) + \theta_2 \sum
g_2(b^*_ib_i) + \theta_1 \sum g_1(b_ib^*_i)\right]$$ 
where the supremum runs
over all pairs of states $(f_1,f_2)$ and $(g_1,g_2)$ on $A$ and $B$
respectively.

Since the right side does not change if we replace $a_i$ by $z_ia_i$ with
$|z_i|=1$, the preceding inequality remains valid with $2\sum |u(a_i,b_i)|$ in
place of $2\left|\sum u(a_ib_i)\right|$ on the left side. Then by a well known
variant of the ``minimax principle'', we find states $f_1,f_2$ and $g_1,g_2$
such that for all $(a,b)$ in $E\times F$ we have
$$2|u(a,b)| \le \alpha_1f_1(a^*a) + \alpha_2f_2(aa^*) + \theta_2g_2(b^*b) +
\theta_1g_1(bb^*)$$
and finally involving (1.3) again we obtain (i).\qed\medskip

\n {\bf 1.4.} For the reader's convenience, we now summarize the basic facts on
the Haagerup tensor product $E\otimes_h F$ of two operator spaces $E,F$. Let
us denote by $E\otimes F$ the algebraic tensor product. Assume $E,F$ are given
together with (completely isometric) embeddings $E\subset A$ and $F\subset B$
in $C^*$-algebras $A$ and $B$. For any $x = \sum a_i\otimes b_i$ in $E\otimes
F$, we define
$$\|x\|_h = \inf\left\{\left\|\sum a_ia^*_i\right\|^{1/2} \left\|\sum
b^*_ib_i\right\|^{1/2}\right\}$$
where the infimum runs over all possible ways to write $x$ as a finite sum
of the form $x = \sum a_i\otimes b_i$. We will denote by $E\otimes_h F$ the
completion of the resulting normed space.

More generally, let $x = [x_{ij}]$ be an $n\times n$ matrix with entries in
$E\otimes F$. We define
$$\|x\|_{(h,n)}  = \inf\{\|y\|_{M_{n,N}(E)} \|z\|_{M_{N,n}(F)}\}$$
where the infimum runs over all $N\ge 1$ and all possible decompositions of $x$
of the form
$$x_{ij} = \sum^N_{k=1} y_{ik} \otimes z_{kj}.$$
By known results (see [CS1-2, PaS]), for a suitable
${\cal H}$, we can find an isometric embedding of
$E\otimes_h F$ into $B({\cal H})$ such that the
above norm $\|~~\|_{(h,n)}$ can be identified
with the norm induced on $M_n(E\otimes F)$ by
$M_n(B({\cal H}))$. This allows us to think of
$E\otimes_h F$ as an operator space. Let $K$ be
any Hilbert space. Then a {\sl bilinear\/} form
$u\colon \ E\times F\to B(K)$ is c.b.\ iff the associated {\sl linear\/} map
$U\colon \ E\otimes_h F\to B(K)$ is well defined and c.b. Moreover we have
$$\|u\|_{cb} = \|U\|_{cb}.\leqno (1.4)$$
It is known ([CS1-2, PaS]) that $\|u\|_{cb}\le 1$
iff there is a Hilbert space $H$ and c.b.\ maps
$\sigma_1\colon \ E\to B(H,K)$ and
$\sigma_2\colon \ F\to B(K,H)$ with
$\|\sigma_1\|_{cb}\le 1$, $\|\sigma_2\|_{cb}\le
1$ such that
$u(a,b) = \sigma_1(a) \sigma_2(b)$ for any $(a,b)$ in $E\times F$.
If $\dim(K) = \infty$, we may actually take $H=K$.

In particular, for a bilinear form $u\colon \ E\times F\to \CC$ we find
$\|u\|_{cb}\le 1$ iff there are $H$ together with $\sigma_1\colon \ E\to
B(H,\CC)$ and $\sigma_2\colon \ F\to B(\CC,H)$ such that $u(a,b) =
\sigma_1(a)\cdot \sigma_2(b)$ for all $(a,b)$ in $E\times F$. We will denote
by $H_c$ and $H_r$ the column and row Hilbertian operator spaces. These are
defined by
$$H_c = B(\CC,H)\quad \hbox{and}\quad H_r = B(H^*,\CC),$$
with the induced operator space structure. Recall (see [ER3]) $(H_c)^* =
(H^*)_r = \ovl H_r$, $(H_r)^* =(H^*)_c = \ovl H_c$.

In particular, in the factorization immediately above we have completely
contractive linear maps $\sigma_1\colon \ E\to \ovl H_r$ and $\sigma_2\colon \
F\to H_c$. Their adjoints $\sigma^*_1 \colon \ H_c\to E^*$ and $\sigma^*_2
\colon \ \ovl H_r\to F^*$ have the same c.b.\ norms.  Thus, the linear map
$\tilde u \colon \ E\to F^*$ associated to $u$ can be factorized through $\ovl
H_r$, i.e.\ we can write $\tilde u = \sigma^*_2\sigma_1$. This leads to a
description of the space $(E\otimes_h F)^*$ in terms of factorization:\ a
bilinear map $u\colon \ E\times F\to \CC$ is c.b.\ iff the linear map $\tilde
u\colon \ E\to F^*$ admits for some $H$ a c.b.\ factorization of the form
$$\tilde u\colon \ E {\buildrel \alpha\over \longrightarrow} H_r {\buildrel
\beta\over \longrightarrow} F^*.$$
Moreover, we have
$$\|u\|_{cb} = \inf\{\|\alpha\|_{cb} \|\beta\|_{cb}\}$$
where the infimum runs over all possible such factorization of $\tilde u$. On
the other hand, for the transposed bilinear form ${}^tu\colon \ F\times E\to
\CC$, we have
$$\|{}^tu\|_{cb} = \inf\{\|\gamma\|_{cb}\|\delta\|_{cb}\}$$
where the infimum runs over all possible c.b.\ factorizations of
$\widetilde{ u}$ of the form
$$\widetilde{ u}\colon \ E {\buildrel \gamma\over\longrightarrow} H_c
{\buildrel \delta\over\longrightarrow} F^*$$
through (this time) a {\sl column\/} space $H_c$.\medskip

\n {\bf 1.5.} We should mention that $E\otimes_h F$ can be realized as a
subspace of the full free product $A*B$ of the $C^*$-algebras, $A,B$
containing $E,F$ respectively, see [CES]. (This is valid also in the unital
case with the unital free product, cf.\ [P2].) In particular, for any
$C^*$-algebra ${\cal B}$ and any finite sum $x = \sum a_i\otimes b_i\otimes
c_i$ in $E\otimes F \otimes {\cal B}$ we have
$$\|x\|_{(E\otimes_h F) \otimes_{\rm min} {\cal B}} = \sup\left\{\left\| \sum
\sigma_1(a_i) \sigma_2(b_i) \otimes c_i\right\|_{\rm min}\right\}$$
where the supremum runs over all $H$ and all pairs of complete contractions
$\sigma_1\colon \ E\to B(H)$, $\sigma_2\colon \ F\to B(H)$.\medskip

\n {\bf 1.6.} The next result from [OP] provides a description of the kind of
bilinear forms that we encounter in this paper.

Let $U \colon \ E\times F\to B(K)$ be a bilinear form. The following are
equivalent:
\item{(i)} There is a decomposition $U = u+v$ with $\|u\|_{cb} +
\|{}^tv\|_{cb} \le 1$.
\item{(ii)} For any $C^*$-algebra ${\cal B}$ and any finite sum $\sum
a_i\otimes b_i \otimes c_i$ in $E\otimes F \otimes {\cal B}$ we have
$$\left\|\sum U(a_i,b_i) \otimes c_i\right\|_{\rm min} \le \sup\left\{\left\|
\sum \sigma_1(a_i) \sigma_2(b_i) \otimes c_i\right\|_{\rm min}\right\}$$
where the supremum runs over all $H$ and all possible pairs $\sigma_1\colon \
E\to B(H)$, $\sigma_2\colon \ F\to B(H)$ of complete contractions {\sl with
commuting\/} ranges.
\item{(iii)} $U$ defines a completely contractive linear map from $(E\otimes_h
F)\cap ({}^tF\otimes_h E)$ equipped with its natural operator space structure.

We will also need the following

\proclaim Proposition 1.7. Let $E,F$ be operator spaces, and let $w\in
E\otimes  F$ (algebraic tensor product). Then there is a finite sequence
$(a_i,b_i)$, $(i=1,2,\ldots, r)$ in $E\times F$ and scalars $\lambda_i>0$ such
that
$$\left\|\sum a_ia^*_i\right\|^{1/2} \left\|\sum b^*_ib_i\right\|^{1/2} =
\|w\|_h\leqno (1.5)$$
and
$$\left\|\sum \lambda_ia^*_ia_i\right\|^{1/2} \left\| \sum \lambda^{-1}_i
b^*_ib_i\right\|^{1/2} = \|{}^tw\|_h. \leqno (1.6)$$

\pf
Without loss of generality we may assume $E$ and $F$ finite dimensional. By
definition of $\|~~\|_h$ (and a compactness argument) there are $(a_i,b_i)$
such that (1.5) holds. Let $\alpha\colon \ E^*\to \ell^r_2$ and $\beta\colon \
\ell^r_2\to F$ be the linear maps defined by
$$\alpha(\xi) = \sum\xi(a_i)e_i\quad \hbox{and}\quad \beta(e_i) = b_i,\leqno
(1.7)$$
so that if $\tilde w \colon \ E^*\to F$ denotes the linear map associated to
$w$, we have
$$w = \beta\alpha.$$
For simplicity we will denote
$$|\alpha|_C = \left\|\sum a^*_ia_i\right\|^{1/2}\quad |\alpha|_R =
\left\|\sum a_ia^*_i\right\|^{1/2}$$
and
$$|\beta|_C = \left\|\sum b^*_ib_i\right\|^{1/2} \quad |\beta|_R = \left\|\sum
b_ib^*_i\right\|^{1/2}.$$
We will use repeatedly the following simple observation:\ for any
$\gamma\colon \ \ell^r_2\to \ell^r_2$ with $\|\gamma\|\le 1$, we have
$$\left\{\eqalign{|\gamma\alpha|_R &\le |\alpha|_R,\quad |\gamma\alpha|_C \le
|\alpha|_C\cr
|\beta\gamma|_R &\le |\beta|_R,\quad |\beta\gamma|_C \le |\beta|_C.}\right.
\leqno (1.8)$$
In particular, this observation allows us to assume that $r$ is the rank of
$w$ and that $(a_i)_{i\le r}$ and $(b_i)_{i\le r}$ are linearly independent.
By the definition of $\|~~\|_h$ again there are maps $\alpha_1\colon \ E^* \to
\ell^r_2$, $\beta_1\colon \ \ell^r_2\to  F$ such that $\tilde w =
\beta_1\alpha_1$ and $|\alpha_1|_C |\beta_1|_R = \|{}^tw\|_h$. By the linear
independence of $(a_i), (b_i)$, there are linear maps $\gamma$ and $\delta$ on
$\ell^r_2$ such that $\alpha_1 = \gamma\alpha$ and $\beta_1 = \beta\delta$.
Moreover since $\beta_1\alpha_1 = \beta\alpha$ we must have $\delta\gamma =
I$, hence $\delta \equiv\gamma^{-1}$. We now write $\gamma$ as a product
$$\gamma = \gamma_1D\gamma_2$$
where $\gamma_1,\gamma_2$ are unitary and $D$ is a diagonal matrix with
coefficients $D_{ii}>0$.

\n  We have then $\delta = \gamma^{-1} = \gamma^{-1}_2
D^{-1} \gamma^{-1}_1$. Hence we can write
$$\tilde w = \beta_1\alpha_1 = (\beta\gamma^{-1}_2D^{-1}) (D\gamma_2\alpha).$$
Now if we replace $\alpha$ and $\beta$ by $\hat\alpha = \gamma_2 \alpha$ and
$\hat\beta = \beta\gamma^{-1}_2$ then (1.8) guarantees that (1.5) still holds,
but on the other hand, setting $\hat\beta_1 = \hat\beta D^{-1}$ and
$\hat\alpha_1 = D\hat\alpha$ we now have
$$|\hat\alpha_1|_C  = |\alpha_1|_C \quad |\hat\beta_1|_R = |\beta_1|_R$$
hence $|\hat\alpha_1|_C |\hat\beta_1|_R = \|{}^tw\|_h$.

\n Moreover, if we denote by $(\hat a_i)$ and $(\hat b_i)$ the sequences
associated to $\hat\alpha$ and $\hat\beta$ as in (1.7) above, then we have
$$\left\|\sum \hat a_i\hat a{}^*_i\right\|^{1/2} \left\|\sum \hat b{}^*_i\hat
b_i\right\|^{1/2} = |\hat \alpha|_R |\hat\beta|_C = \|w\|_h$$
and on the other hand
$$\left\|\sum D^2_{ii}\hat a{}^*_i \hat a_i\right\|^{1/2} \left\|\sum
D^{-2}_{ii} \hat b_i\hat b{}^*_i\right\|^{1/2} = |\hat\alpha_1|_C
|\hat\beta_1|_R = \|{}^tw\|_h.$$
Hence letting $\lambda_i= D_{ii}$ we obtain the announced result.\qed

\ms
\n {\bf 1.8.} A $C^*$-algebra $A$ is called WEP (for weak expectation
property) if the inclusion map $A\to A^{**}$ can be factorized completely
contractively through $B(H)$ for some $H$. This notion goes back to Lance [L].
Many equivalent definitions are known. The one we will use is the following
characterization due to Kirchberg [Ki2]:\ let $C$ be the (full) $C^*$-algebra of the
free group with (say) countably infinitely many generators. Then $A$ is WEP
iff $A \otimes_{\rm min} C = A \otimes_{\rm max} C$ (that is to say the
minimal and maximal $C^*$-norms coincide on the algebraic tensor product 
$A\otimes C$). A simpler proof was
given in [P2].

Following Kirchberg [Ki1], we will say that a $C^*$-algebra $A$ is QWEP if it
is a quotient of a WEP $C^*$-algebra. It is an outstanding open question
 whether {\sl every\/} $C^*$-algebra is QWEP. This is equivalent
to Connes' question whether every von Neumann algebra (on a separable
Hilbert space) embeds into an ultra-product of the hyperfinite factor
 (see [Ki1]
for more on this). 
\ms
\n {\bf 1.9.} Let $E$ be a finite dimensional operator space. Let $C  =
ex(E)$. Then for any $C^*$-algebra $B$ and any (closed 2-sided) ideal 
$I\subset B$, we have a canonical isomorphism
$$T\colon \ E\otimes_{\rm min} (B/I) \to (E\otimes_{\rm min}
B)/{(E\otimes_{\rm min} I)}$$
with $\|T\| \le C$ (and obviously $\|T^{-1}\|\le 1$). 

\n A $C^*$-algebra is exact iff  $ex(E)=1$ for any finite dimensional
subspace $E \subset A$. We will use this in the following manner: let
$I$ and $B$ be as above. Fix $\vp>0$.  Then, for any  $t \in
A\otimes (B/I)$ (algebraic tensor product)
there is a lifting $\hat t \in A\otimes B$ (again {\it algebraic} tensor product)
with $\|\hat t\|_{\min} \le (1+\vp) \|t\|_{\min}.$  See [P2] (or [P1, ER1])
for details.

\n { \bf Remark 1.10.} We draw the reader's attention to the fact that, if
$A$ is arbitrary, since $A\otimes_{\min} (B/I)$ is a quotient 
$C^*$-algebra of $A\otimes_{\min} B$, there always exist a lifting
of $t $ 
(say $\hat t$ with $\|\hat t\|_{\min} \le (1+\vp) \|t\|_{\min}$)
in the { \it completed}  tensor product $A\otimes_{\min} B$.
 But this seems to be of little use for us in the sequel.
Indeed, it is always true that $t$ admits a lifting  
$\hat t \in A\otimes_{\max} B$ with $\|\hat t\|_{\max} \le (1+\vp) \|t\|_{\min}$,
but in general
we cannot derive from this
that  $\|  t\|_{\max} \le (1+\vp) \|t\|_{\min}$,
unless we know that $\hat t$ is in the algebraic tensor
product, in which case we {\it do} obtain
 $\|  t\|_{\max} \le (1+\vp) \|t\|_{\min}$.
The crucial point is that, in general, 
the canonical map  $A\otimes_{\max} B\to A\otimes_{\min} B$ is 
{\it not injective}   unless restricted to the {\it algebraic} tensor product.

\medskip

\n {\bf \S 2. Proofs}

\n We will use exactness in the same way as in [P3] and [JP], through the
following result which is implicit in [JP].

\proclaim Lemma 2.1. Let $E,F$ be exact operator spaces and let $C = ex(E)
ex(F)$. Let $A_1,A_2$ be two $C^*$-algebras such that either $A_1$ or $A_2$ is
QWEP. Let $u\colon\ E\times F\to B({\cal H})$ be a j.c.b.\ bilinear map. Then
for any finite sequences $(a_i)$, $(b_j)$, $(x_i)$, $(y_j)$ in $E,F,A_1,A_2$
respectively we have
$$\left\|\sum u(a_i,b_j) \otimes x_i\otimes y_j\right\|_{ B({\cal
H}) \otimes_{\rm min} (A_1 \otimes_{\rm max} A_2)} \le C\left\|\sum a_i
\otimes x_i\right\|_{E\otimes_{\rm min}A_1} \left\|\sum b_j\otimes
y_j\right\|_{F\otimes_{\rm min} A_2}.\leqno (2.1)$$
Equivalently, $u$ ``extends'' to a bounded bilinear form from $E\otimes_{\rm
min} A_1\times F \otimes_{\rm min} A_2$ to $B({\cal H}) \otimes_{\rm min} (A_1
\otimes_{\rm max} A_2)$, taking $(e\otimes a_1, f\otimes a_2)$ to $u(e,f)
\otimes a_1\otimes a_2$. 

\pf We may assume that $A_1$ is QWEP, so that $A_1\simeq B_1/I_1$ with
$B_1$ WEP. Taking for $B_2$ the full $C^*$-algebra of a suitably large
free group we may assume $A_2 \simeq B_2/I_2$. We denote by $q_i\colon \
B_i\to A_i$ the quotient map. We will use the isometric identity 
$$B_1\otimes_{\rm min} B_2 \simeq B_1 \otimes_{\rm max} B_2.$$
This is due to Kirchberg [Ki2], (see [P2] for a simpler proof). In particular
we have
$$\|q_1\otimes q_2\colon \ B_1 \otimes_{\rm min} B_2 \to A_1 \otimes_{\rm max}
A_2\|_{cb} \le 1,\leqno (2.2)$$
hence
$$\|I_{B{\cal H}) } \otimes q_1\otimes q_2\colon \ {B({\cal H})
\otimes_{\rm min} (B_1\otimes_{\rm min} B_2)} \to {B({\cal H})
\otimes_{\rm min} (A_1\otimes_{\rm max} A_2)}\|  \le 1,\leqno (2.2)'$$

We may as well assume, without loss of generality, that $E$ and $F$ are finite dimensional. 
We now use the exactness of $E$ and $F$ (see (1.9)):\ let $C_1 = ex(E)$ and $C_2 = ex(F)$, assuming
$$\left\|\sum a_i\otimes x_i\right\|_{\rm min} < 1\quad \hbox{and}\quad
\left\|\sum b_j \otimes y_j\right\|_{\rm min} < 1$$
we can find elements $t_1 \in E\otimes B_1$ and $t_2\in F\otimes B_2$ with
$\|t_1\|_{\rm min} < C_1$, $\|t_2\|_{\rm min} < C_2$ such that $(I\otimes
q_1)(t_1) = \sum a_i\otimes x_i$ and $(I\otimes q_2)(t_2) = \sum b_j \otimes
x_j$. By (1.1)$'$   we have
$$\|(u)_{B_1,B_2} (t_1,t_2)\| \le \|u\|_{jcb} \|t_1\|_{\rm min} \|t_2\|_{\rm
min},$$
and moreover $(u)_{B_1,B_2} (t_1,t_2)$ lies in the algebraic
tensor product ${B({\cal H})
\otimes  B_1\otimes B_2}$.  
But clearly
$$(I\otimes q_1\otimes q_2)(u)_{B_1,B_2} = (u)_{A_1,A_2} \circ (I\otimes q_1,
I\otimes q_2)$$
hence by (2.2)' (and Remark 1.10)
$$\left\|\sum u(a_i,b_j) \otimes x_i\otimes y_j\right\|_{B({\cal H})
\otimes_{\rm min} (A_1\otimes_{\rm max} A_2)} \le C_1C_2\|u\|_{jcb}.$$
By homogeneity, this completes the proof.\qed

\proclaim Lemma 2.2. Let $A,B$ be {\sl arbitrary\/} $C^*$-algebras and let
$u\colon \ A\times B\to \CC$ be a j.c.b.\ bilinear map of finite
rank with $\|u\|_{jcb}\le 1$. Then $u$ satisfies the conclusion of
Lemma~2.1 with $C=1$.

\pf The same argument as above for Lemma~2.1 can be used to establish this,
but the exactness of $E$ and $F$ has to be replaced by the following known
fact:\ Let $E\subset A$ be a finite dimensional subspace and let $Z$ be a
finite dimensional operator space. Then for any linear map $\alpha\colon \
A\to Z$ with $\|\alpha\|_{cb}<1$, the restriction $\alpha_{|E}\colon \ E\to Z$
is 1-exact in the following sense:\ there is an integer $N$, a subspace
$G\subset M_N$ and a factorization $E {\buildrel v\over
\longrightarrow} G {\buildrel w\over\longrightarrow} Z$ of $\alpha_{|E}$ with
$\|v\|_{cb} \|w\|_{cb} < 1$.

\n This was observed in [P1]. Here is a quick sketch
of proof:\ let $I\subset B$ be an ideal in a
$C^*$-algebra $B$ and let $q\colon \ B\to B/I$ be
the quotient map. For any operator space $Z$ we
denote by $q_Z\colon \ Z\otimes_{\rm min} B\to
Z\otimes_{\rm min} B/I$ the (contractive) map
associated to $I_Z \otimes q$. Moreover we let
$$Q(Z) = (Z \otimes_{\rm min} B)/(Z \otimes_{\rm min} I)\quad \hbox{ and }\quad
R(Z) = Z \otimes_{\rm min}(B/I).$$
It suffices to show (and this is precisely what we use to prove Lemma~2.2)
that $a_{|E}\colon \ E\to Z$ naturally induces a mapping $\hat \alpha\colon \
R(E)\to Q(Z)$ with $\|\hat \alpha\| < 1$. To verify this, note that since
$q_A$ is a $*$-homomorphism, it induces an isometric isomorphism
$$R(A) \simeq A \otimes_{\rm min} B/\ker(q_A).\leqno (2.3)$$
But since $\alpha$ has finite rank,
 it is easy to see that   the map $\alpha\otimes
q\colon \ A 
\otimes_{\rm min} B\to Z \otimes_{\rm min} (B/I)$ (obtained by extension from
$\alpha\otimes q$) must vanish on $\ker(q_A)$, hence by (2.3) it defines a map
$\tilde\alpha\colon \ R(A)\to Q(Z)$ with $\|\tilde\alpha\|\le \|\alpha \otimes
q\| \le \|\alpha\|_{cb} < 1$.
Finally since $\hat\alpha\colon \ R(E)\to Q(Z)$
is but the restriction of
$\tilde\alpha$ to $R(E)$ we obtain $\|\hat\alpha\| \le \|\alpha\|_{cb}<1$.

\n  Using this fact we complete the proof as follows.
Assume for simplicity that $\|u\|_{jcb}<1$. Since $u$ is assumed of finite
rank, there is a finite dimensional operator space $Z$ and linear maps
$\alpha\colon \ A\to Z$ and $\beta\colon \ B\to Z^*$ with $\|\alpha\|_{cb} \le
1$ and $\|\beta\|_{cb} \le 1$, such that
$$u(a,b) = \langle\alpha(a), \beta(b)\rangle$$
Let $E\subset A$ and $F\subset B$ be finite dimensional subspaces (we can take
$E = \hbox{span}(a_i)$ and $F = \hbox{span}(b_j))$. Then the maps $\alpha_{|E}
\colon \ E\to Z$ and $\beta_{|F}\colon \ F\to Z^*$ are 1-exact in the above
sense. The rest of the proof is then identical to that of Lemma~2.1.\qed

\proclaim Proposition 2.3. Let $A$ be a $C^*$-algebra and let $F$ be an
exact  operator space. Then any j.c.b. bilinear form
$u\colon \ A\times F\to \CC$   with $\|u\|_{jcb}\le 1$ satisfies
the conclusion of
Lemma~2.1 with $C=ex(F)$. 

\pf Indeed, we can run the same argument
as for Lemma 2.1, but using the remarks at the end of \S 1.9.
With the same notation as in Lemma 2.1 (but here $E=A$
and $C_1=1$)
we can find   $t_1$ in the completed tensor product $  A\otimes_{\min} B_1$
satisfying all the properties
in the proof of Lemma 2.1, but since $t_2$ can still be chosen
in the algebraic tensor product $F\otimes B_2$, 
we find again 
$ (u)_{B_1,B_2} (t_1,t_2)$ in the algebraic tensor product $ B_1 \otimes B_2$,
hence we can complete the proof
as for Lemma 2.1. \qed

\medskip

In the previous versions of GT, the crucial ingredient is always the existence
of a ``special'' realization of Hilbert space as a function space of some sort.
The analog of this in our situation is the span of a free family of
``generalized circular elements'' in the sense of [S1]. This is an
extension of Voiculescu's circular systems (see [VDN]), as
follows.

Let $H$ be any Hilbert space. Let $H^{\otimes n} = H\otimes_2\cdots \otimes_2
H$ ($n$ times). We denote by ${\cal F}(H)$ the full Fock space over $H$, i.e.
$${\cal F}(H)= \CC \oplus H \oplus H^{\otimes 2} \oplus\cdots\oplus H^{\otimes
n}\oplus\cdots$$
We denote by $\Omega$ the unit in $\CC$ view as an element of ${\cal F}(H)$.
For any $h$ in $H$, we denote by 
$\ell(h)\colon \ {\cal F}(H)\to {\cal F}(H)$ (resp.\ $r(h)\colon {\cal
F}(H)\to{\cal F}(H)$) the left (resp.\ right) creation operator, defined by:\
$\ell(h)\Omega = h$ (resp.\ $r(h)\Omega = h$) and for any $x$ in $H^{\otimes
n}$ with $n>0$:\ $\ell(h)x = h\otimes x$ (resp.\ $r(h) x = x\otimes h$). We
will assume that $H$ admits an orthonormal basis which can be split in two
parts with equal cardinality, 
$$\{e_i\mid i\in I\}\quad\hbox{and}\quad \{e'_i\mid i\in I\}$$
so that the union $\{e_i\mid i\in I\} \cup \{e'_i\mid i\in I\}$ is an
orthonormal system. We will denote
$$\eqalign{\ell_i = \ell(e_i) &\quad \hbox{and}\quad r_i = r(e_i)\cr
\ell'_i = \ell(e'_i) &\quad \hbox{and}\quad r'_i = r(e'_i).}$$
Then we define, for any $\lambda>0$, the ``generalized circular elements'' as
follows:
$$\eqalignno{
c_i(\lambda) &= \lambda^{1/2}\ell_i + \lambda^{-1/2}\ell'{}^*_i\cr
\noalign{\hbox{and}}
d_i(\lambda) &= \lambda^{1/2} r'_i + \lambda^{-1/2} r^*_i.}$$

The von Neumann algebra generated by such  systems is studied in [S1].
 (When $\lambda=1$, we recover circular elements in Voiculescu's sense.)

We will use the following basic properties of these operators
 (the key point is the second one):

\proclaim Lemma 2.4. Let $\lambda_i>0$ be fixed $(i\in I)$. Let $x_i =
c_i(\lambda_i)$ and $y_i = d_i(\lambda_i)$.
\item{(i)} Let $E,F$ be operator spaces and let $(a_i,b_i)_{i\in I}$ be a
finitely supported family in $E\times F$. Then we have 
$$\eqalign{\left\|\sum a_i\otimes x_i\right\| &\le \left\|\sum \lambda_i
a^*_ia_i\right\|^{1/2} + \left\|\sum \lambda^{-1}_ia_ia^*_i\right\|^{1/2}\cr
\left\|\sum b_i\otimes y_i\right\| &\le \left\|\sum
\lambda_ib^*_ib_i\right\|^{1/2} + \left\|\sum\lambda^{-1}_i
b_ib^*_i\right\|^{1/2}.}$$
\item{(ii)} The families $\{x_i\} \{y_j\}$ ``double commute'' which means:\
$x_iy_j = y_jx_i$ and $x^*_iy_j = y_jx^*_i$ for all $i,j$ in $I$.
\item{(iii)} $\langle x_iy_i\Omega,\Omega\rangle = 1$ and $\langle x_iy_j
\Omega,\Omega\rangle = 0$ for all $i\ne j$.

\pf To prove (i) note that each $\{\ell_i\}$, $\{\ell'_i\}$, $\{r_i\}$ or
$\{r'_i\}$ is a family of isometries with orthogonal ranges and if $(s_i)$ is
any such family we must have $\left\|\sum s_is^*_i\right\| \le 1$. By an easy
and well known consequence of Cauchy-Schwarz, we have for any finitely
supported family of operators $(s_i,t_i)_{i\in I}$ 
$$\eqalignno{
\left\|\sum s_i\otimes t_i\right\| &\le \left\|\sum s_is^*_i\right\|^{1/2}
\left\|\sum t^*_it_i\right\|^{1/2}\cr
\noalign{\hbox{and of course also}}
&\le \left\|\sum s^*_is_i\right\|^{1/2} \left\|\sum t_it^*_i\right\|^{1/2}.}$$
Thus by the triangle inequality we have
$$\eqalign{\left\|\sum a_i\otimes x_i\right\|&\le \left\|\sum
a_i\lambda^{1/2}_i \otimes\ell_i\right\| + \left\|\sum a_i\lambda^{-1/2}_i
\otimes \ell'{}^*_i\right\|\cr
&\le \left\|\sum a^*_ia_i\lambda_i\right\|^{1/2} \left\|\sum \ell_i\ell^*_i
\right\|^{1/2} + \left\|\sum a_ia^*_i\lambda^{-1}_i\right\|^{1/2} \left\| \sum
\ell'_i\ell'{}^*_i\right\|^{1/2}}$$
which implies the first part of (i). The same argument proves the other part. 

\n To prove (ii) note that for all $h,k$ in $H$ the operators $\ell(h)$ and
$r(k)$ obviously commute and 
$$r(k)^* \ell(h) - \ell(h) r(k)^* = \langle h, k\rangle   P_\Omega.$$
It is then but an elementary verification to check (ii).

\n (iii)\ We have $y_i\Omega = \lambda^{1/2}_i e'_i$ and $x_ie'_i =
\lambda^{1/2}_i e_i\otimes e'_i + \lambda^{-1/2}_i\Omega$ hence $\langle
x_iy_i\Omega,\Omega\rangle = \langle\Omega,\Omega\rangle = 1$. Similarly if
$i\ne j$ we find $x_iy_j\Omega = \lambda^{1/2}_i\lambda^{1/2}_j e_i\otimes
e'_j$ hence $\langle x_iy_j\Omega,\Omega\rangle = 0$.\qed

 In order to be able to
use the generalized circular families, it is crucial to know the following:

\proclaim Lemma 2.5. The von Neumann algebra $W^*(x_i : i\in I)$ 
generated by the family $\{x_i
\mid i\in I\}$ defined in the preceding lemma is QWEP.

\pf Let $H_{\RR}$ be a real Hilbert space of dimension $2|I|$, and denote by ${f_i, f'_i}$,
$i\in I$ its orthonormal basis. Let $U_t : H_{\RR}\to H_{\RR}$ be given by
$U_t(f_i) = \cos(\log(\lambda^2_i) t)f_i + \sin(\log(\lambda^2_i) t)f'_i$,
$U_t (f'_i) = \cos(\log(\lambda^2_i t))f'_i - \sin (\log(\lambda^2_i) f_i$.  Then $U_t$ is a
one-parameter group of orthogonal transformations on $H_{\RR}$.  Let
$M= \Gamma(H_{\RR},U_t)''$ in its GNS representation with
respect to the free quasi-free state 
$\phi_U$ (see Definition 2.3 of
[S1]). By [S1, Section
4], we find that $M$ can be viewed as generated
by the operators $l(e_i) + \lambda_i l(e_i')^*$
and is therefore isomorphic to the algebra
$W^*(x_i : i\in I)$.  Hence it is sufficient to prove that $M$ is QWEP.

Since $U_t$ is almost-periodic, it follows from Theorem 6.4 and Theorem 6.7 of
[S1] that  $$(M,\phi_U)\cong *_{i\in I} \left\{
\left( M_{2\times 2}, \hbox{Tr}\left(\left(
 \matrix{  {{1}\over{1+\lambda^2}} & 0 \cr 0 &
 {{\lambda^2}\over{1+\lambda^2}} }
\right) \cdot \right)\right) * L^\infty[0,1]
\right\}.$$  
Hence by the results of K. Dykema
[D], the centralizer
$M^{\phi_U}$ of $\phi_U$ is isomorphic to $L(\FF_\infty)$.  Thus $M^{\phi_U}$ is a
factor.   Moreover, by Kirchberg's results,
it is QWEP, since $L(\FF_\infty)$
 can be embedded into the
ultrapower of the hyperfinite II$_1$
 factor (see [Ki1]).

Since $U_t$ is almost-periodic, it follows that the modular group of $\phi_U$ is almost-periodic,
so that $\phi_U$ is an almost-periodic state
 (see [C]). By Connes'
results in [C], it follows from the fact that
$M^{\phi_U}$ is a factor that
$$ M \cong \left( M^{\phi_U} \otimes B(\ell^2)  \otimes B(H)\right)
\rtimes G,$$ where
$G$ is a discrete group (isomorphic to the multiplicative subgroup of $(0,+\infty)$ generated
by the set $\{\lambda^2_i :i \in I\}$). Since $G$ is Abelian,
hence amenable, and $M^{\phi_U}$ is QWEP,  it follows that  $M$ is QWEP (see the remark
after Prop. 1.3 in [Ki1]). \qed

\n { \bf Remark.} It is possible to give an alternate proof that $M$ is
QWEP. Indeed,
it is sufficient to find a sequence of states $\phi_i$ on matrix
algebras $A_i = M_{n_i
\times n_i}$, a free ultrafilter $\omega$ and an embedding
$$i : M \to \prod_\omega (A_i, \phi_i) $$ so that $\prod \sigma^{\phi_i}
|_{i(M)} = \sigma^{\phi_U}$.
Indeed, the latter condition implies that there exists a
state-preserving conditional
expectation $E: \prod_\omega(A_i,\phi_i)$ onto $i(M)$; since $\prod A_i$
is WEP, $\prod_\omega
(A_i,\phi_i)$ is QWEP, and thus $i(M)$ is QWEP.

\n To construct the embedding $i$ one can utilize the model for free
quasi-free
states involving matrices with CAR variables as entries (described in
[S2]).
To make sure that the matrices involved stay bounded in norm, one must
cut them off using continuous functional calculus; it can be verified that this procedure
can be performed
in a way that does not affect their joint $*$-distribution.  We leave
the details
to the reader.

We can now prove our main results.

\n {\bf Proof of Theorem 0.3.} Let $A_1$ (resp.\ $A_2$) be the von~Neumann
algebra generated by $(x_i)$ (resp.\ $(y_i)$). By Lemma~2.4 (iii)
$$\eqalignno{\left|\sum U(a_i,b_i)\right| &= \left|\sum_{ij} U(a_i,b_j)
\langle x_iy_j\Omega,\Omega\rangle\right|\cr
\noalign{\hbox{hence by Lemma 2.4 (ii)}}
&\le \left\|\sum_{ij} U  (a_i,b_j) \otimes x_i\otimes y_j\right\|_{B({\cal H})
\otimes_{\rm min} (A_1 \otimes_{\rm max} A_2)}\cr
\noalign{\hbox{hence by Lemma 2.1}}
&\le C\left\|\sum a_i \otimes x_i\right\|_{\rm min} \left\|\sum b_i\otimes y_i
\right\|_{\rm min}}$$
hence by Lemma 2.4 (i) we obtain (0.13). Taking $\lambda_i=\lambda$ for all
$i$ and then choosing $\lambda^2  = \left\|\sum b_ib^*_i\right\|^{1/2}
\left\|\sum b^*_ib_i\right\|^{-1/2}$, we derive (0.14) from (0.13).\qed

\n {\bf Proof of Theorem 0.4.} By Theorem 0.3 and by (1.3) we have (note that
the argument of $U(a_i,b_i)$ can be absorbed e.g.\ by $a_i$)
$$\sum|U(a_i,b_i)| \le C[I+II]$$
where
$$\eqalign{I &= \left\|\sum \lambda_ia^*_ia_i\right\| + \left\|\sum
\lambda^{-1}_i b_ib^*_i\right\|\cr
II &= \left\|\sum \lambda^{-1}_i a_ia^*_i\right\| + \left\|\sum \lambda_i
b^*_ib_i\right\|.}$$
By a very easy adaptation of Proposition~1.3 this implies the existence of
states $f_1,f_2,g_1,g_2$ on $A$ and $B$ such that for any $(a,b)$ in $E\times
F$ and any $\lambda>0$ (note it is crucial that the states do not depend on
$\lambda$!) 
$$|U(a,b)| \le C[\lambda f_1(a^*a) + \lambda^{-1}g_1(bb^*) + \lambda^{-1}f_2(aa^*) +
\lambda g_2(b^*b)].$$
Applying this to $(ta,t^{-1}b)$ instead of $(a,b)$,  and taking the infimum over $t>0$ we obtain by
(1.3)
$$\eqalign{|U(a,b)| &\le 2C(\lambda f_1(a^*a) + \lambda^{-1} f_2(aa^*))^{1/2}
(\lambda^{-1}g_1(bb^*) + \lambda g_2(b^*b))^{1/2}\cr
&\le 2C(f_1(a^*a) g_1(bb^*) + f_2(aa^*)g_2(b^*b) + R_\lambda)^{1/2},}$$
where
$$R_\lambda = \lambda^2f_1(a^*a) g_2(b^*b) + \lambda^{-2}f_2(aa^*) g_1(bb^*).$$
By (1.3) we have
$$\eqalignno{\inf_{\lambda>0} R_\lambda &= 2(f_1(a^*a) g_2(b^*b) f_2(aa^*)
g_1(bb^*))^{1/2}\cr
\noalign{\hbox{hence by (1.3) again (with a different grouping of terms}}
&\le f_1(a^*a) g_1(bb^*) + f_2(aa^*) g_2(b^*b)\cr
\noalign{\hbox{thus we  obtain}}
|U(a,b)| &\le 2\sqrt 2\ C(f_1(a^*a) g_1(bb^*) + f_2(aa^*) g_2(b^*b))^{1/2}}$$
which implies (0.15).
Clearly, (0.15) implies (0.16) by Cauchy-Schwarz.

\n Now, by Proposition 1.7, (0.16) can be
reformulated as follows:\ for any $w$ in $E\otimes F$ we have
$$|\langle U,w\rangle| \le K [\|w\|_h + \|{}^tw\|_h].$$
Thus $U$ defines a continuous linear form with norm $\le K$ on the
subspace
$$\{(w, {}^tw)\mid w\in E\otimes F\} \subset (E\otimes_h F) \oplus_1 (F
\otimes_h E)$$
equipped with the norm $\|(x,y)\| = \|x\|_h + \|y\|_h$. By the Hahn-Banach
theorem, there are linear forms $\varphi_1\in (E\otimes_h F)^*$ and
$\varphi_2\in (F\otimes_h E)^*$ with $\max\{\|\varphi_1\|, \|\varphi_2\|\} \le
K$ such that
$$\langle U,w\rangle = \varphi_1(w) + \varphi_2({}^tw).\leqno \forall~w\in
E\otimes F$$
Let $u$ (resp. ${}^t v$) be the bilinear
 forms associated to $\varphi_1$ (resp. $\varphi_2$). Going back to (0.5),
we  see that $\|\varphi_1\|\le K$
(resp.  $\|\varphi_2\|\le K$ ) is equivalent
to
$\|u\|_{cb}\le K$, (resp.  $\|{}^tv\|_{cb}\le K$). Clearly $U=u+v$. 
Thus the final assertion follows.
\qed

\n {\bf Proof of Theorem 0.5.} This clearly reduces to the finite rank case.
The proof is then the same as for Theorem~0.4 but using Lemma~2.2 instead of
Lemma~2.1.\qed

\n {\bf Proof of Corollary 0.6.} Let $\varphi_1,\varphi_2$ be as in the  above
 proof of Theorem 0.4.
Now since $\otimes_h$ is an injective tensor product (\cf [ER1]), $\varphi_1$ and
$\varphi_2$ can be extended to linear forms 
$$\Phi_1\colon\ A\otimes_h B \to \CC\quad \hbox{and}\quad \Phi_2\colon\ B\otimes_h
A\to \CC$$
such that $\|\Phi_i\| = \|\varphi_i\|$, $(i=1,2)$. Let $U_i\colon \ A\times B\to
\CC$ be the bilinear forms associated to $\Phi_i$ \  $(i=1,2)$. By (0.8) and (1.4),
we have
$$\eqalign{\|U_1 + U_2\|_{jcb} &\le \|\Phi_1\| + \|\Phi_2\|\cr
&\le 4\sqrt 2\ C.}$$
Thus if we set
$\hat U = U_1 + U_2$, we obtain Corollary 0.6.\qed\medskip

\n {\bf Proof of Corollary 0.7.} 
Let $u$ be the bilinear form
associated to $T$, so that $T=\tilde u$. By the proof of Corollary 0.6, we have a
decomposition
$$u  = u_1+u_2$$
with 
$$\max(\|u_1\|_{cb}, \|{}^tu_2\|_{cb}) \le 2\sqrt 2\ C \|u\|_{jcb}.$$
By the results recalled in 1.4, the linear map $\tilde u_1\colon \ E\to F^*$
(resp.\ $\tilde u_2\colon \ E\to F^*$) can be factorized through $H_r$ (resp.\
$H_c$) for some $H$. More precisely, we have factorizations
$$\eqalign{\tilde u_1\colon \ &E {\buildrel \alpha_1 \over \longrightarrow} H_r
{\buildrel \beta_1 \over \longrightarrow} F^*\cr
\tilde u_2\colon \ &E {\buildrel \alpha_2 \over \longrightarrow} H_c
{\buildrel \beta_2 \over \longrightarrow} F^*}$$
with
$$\max\{\|\alpha_1\|_{cb} \|\beta_1\|_{cb}, \|\alpha_2\|_{cb}
\|\beta_2\|_{cb}\} \le 2\sqrt 2\ C \|u\|_{jcb.}$$
By homogeneity, we can adjust $\alpha_i, \beta_i$ so that
$$\|\alpha_1\|_{cb} = \|\beta_1\|_{cb}\quad \hbox{and}\quad \|\alpha_2\|_{cb}
= \|\beta_2\|_{cb}.$$
Then, if we define
$$v\colon \ E\to H_r\oplus H_c \quad  \hbox{and}\quad w\colon \ H_r \oplus
H_c\to F^*$$
by
$$v(e) = \alpha_1(e) \oplus \alpha_2(e)\quad \hbox{and}\quad w(x\oplus y) =
\beta_1(x) + \beta_2(y),$$
  we obtain $\tilde u = wv$ and
$$\eqalign{\|w\|_{cb} \|v\|_{cb} &\le \max\{\|\alpha_1\|_{cb},
\|\alpha_2\|_{cb}\} \max\{\|\beta_1\|_{cb}, \|\beta_2\|_{cb}\}\cr
&\le 2\sqrt 2\ C \|u\|_{jcb}.}$$
\qed

\n {\bf Proof of Corollary 0.8.} We prove this
using duality. By (0.16) and by Proposition~1.7
  for any $w$ in $E\otimes F$ and any $U$
as in (0.16) we have
$$\eqalign{|\langle w,U\rangle| &\le
2^{3/2}C[\|w\|_h + \|{}^tw\|_h]\cr &\le 2^{5/2}
C\max\{\|w\|_h, \|{}^tw\|_h\}}$$ hence by (0.11)
$$\|w\|_\wedge = \sup\{|\langle w,U\rangle|\
\|U\|_{jcb}\le 1\} \le 2^{5/2} C
\max\{\|w\|_h, \|{}^tw\|_h\},$$
which establishes the first point.
The same argument allows us to deduce the second point from Theorem~0.5.
Finally, the  third point  also follows from Theorem~0.5 using
the fact that $A^* \otimes_{\rm min} B^*$  can  be identifed with the closure
of the finite rank maps in  $CB(A,B^*)$.\qed

\n {\bf Remark 2.6.} By Remark 2.3, if either $A$ or $B$ is exact,
the conclusions of  Theorem 0.5 and Corollary 0.7 are valid
without any approximability assumption. We do not know whether
they are valid in full generality. Note however  that, by our results,
an operator $T \colon\ A \to B^*$ factors   through a space 
  of the form $  H_r\oplus K_c$   iff it is approximable pointwise by a
  net of finite rank maps unifomly bounded in $CB(A,B^*)$. Indeed, since the identity
of $  H_r\oplus K_c$ is obviously approximable in this way,   the only if part also holds.
Thus,  if Theorem 0.5 holds without  any approximability assumption, this means
that any $T$ in $CB(A,B^*)$ is approximable pointwise by a
  unifomly   bounded net of finite rank maps.

\medskip
\n {\bf \S 3. Some applications}

\n In classical Banach space theory, $GT$ has several well known consequences
(see e.g.\ [P5]). Our result ``automatically'' allows to transfer some of these
to the operator space setting. For instance, let  $E$ be a Banach space. It is known
that     $E$  and its dual $E^*$  both embed in  an $L_1$-space iff $E$ is
isomorphic to a Hilbert space. Actually, this remains  true
if we replace ``$L_1$-space" by ``non-commutative $L_1$-space"
(see [P5]).  In the operator space case, we could not obtain
an analogous characterization, but the next result essentially 
reduces the problem to the class of subspaces of quotients
of $H_r\oplus K_c$. Recall that an operator space
$E$ is said to have the CBAP (for completely bounded aproximation
property)  if the identity of $E$ is the pointwise limit of
a net of finite rank maps, with uniformly bounded c.b. norms.   

\proclaim Corollary 3.1. Let $E$ be an operator space such that:
\item{(i)} $E$ and $E^*$ both embed completely isomorphically in a
non-commutative $L_1$-space (meaning the predual of a von~Neumann algebra),
\item{(ii)} $E$ has the $CBAP$,\hfil\break
   then $E$ is completely isomorphic to a quotient of a subspace of
$H_r\oplus K_c$ for some Hilbert spaces $H,K$.

\pf Note that by the  Banach space result just recalled, $E$ must be
isomorphic to a Hilbert space (in particular it is reflexive).
Let $A,B$ be $C^*$-algebras. Let $J_1\colon \ E^*\to A^*$ and $J_2\colon \
E\to B^*$ be completely isomorphic embeddings. Let $Q_1\colon \ A\to E$ be
the adjoint of $J_1$. Then let $T = J_2 Q_1\colon \ A\to B^*$. By Corollary
0.7, $T$ factors through a space $X$ of the form $X = H_r\oplus K_c$ for
suitable $H,K$. Going back to $E$, we easily deduce from this that $E$ is
completely isomorphic to a quotient of a subspace of $X$.\qed

\n {\bf Remark.} In general, we do not know whether the assumption (ii) is needed.
Note however that if either $E$ or $E^*$ is exact, then 
using Remark 2.3, (ii)  can be dispensed with.

\n {\bf Remark.} Let $E$ be the closed span of the classical Rademacher
functions in $L_1([0,1], dt)$. By  a result
due to Lust-Piquard and the first author (see [P6, p. 107]),  the dual $E^*$ is
completely isomorphic to a subspace of $R\oplus C$, namely to the closed span in
$R\oplus C$ of $\{e_{1i} \oplus e_{i1} \mid i\ge 1\}$. Let us denote by $S_1$ the predual of
$B(\ell_2)$. Note that
$R$ and $C$ both embed in $S_1$ ($R$ can be identified naturally with the columns
in $S_1$, and $C$ with the rows, in a suitable duality between $S_1$  and
$B(\ell_2)$.)
Thus, $E$ is a good illustration of the preceding corollary. In this case we have
$E\subset L_1([0,1], dt)$ and $E^*\subset S_1\oplus S_1$ where $S_1$ is the
predual of $B(\ell_2)$.

\n Another illustration is provided more generally by the family  $\{x_i\mid i\in
I\}$ appearing in Lemma~2.4: assuming for simplicity that 
$\lambda_i=\lambda>0$ does not depend on $i$,
let  $E_\lambda$ denote  the
weak-$*$ closure of $\{x_i\mid i\in I\}$ in $B({\cal F}(H))$.
Note that again  $E_\lambda$ embeds into $R\oplus C\subset  S_1\oplus S_1$ 
(by Lemma 2.4 (i)) and moreover one can show 
there is a completely bounded projection 
from $W^*(x_i : i\in I)$  onto $E_\lambda$. Hence
${E_\lambda}^*$  also embeds in a non-commutative $L_1$-space.

\n {\bf Remark.} The preceding topic is also of interest
in the isometric case. Indeed, by [Sc]  there are 
examples of finite dimensional   Banach spaces $E$,
   not isometric to   Hilbert spaces, but   such that
both $E$ and $E^*$ embed isometrically into $L_1$. Curiously, however,
this is   known {\it in the real case only}. The complex case apparently
remains open, as well as the infinite dimensional  one (either real or
complex).  Analogously, we do not have any satisfactory completely
isometric version of Corollary 3.1. 
 
The next result gives a characterization of the operator spaces $E$ such that 
both $E$ and $E^*$ are exact
 (see [OP] for a different characterization of the same  class of  spaces).   

\proclaim Corollary 3.2. An operator space $E$ is exact as well as its
(operator space) dual $E^*$ iff $E$ is completely isomorphic to  
 $H_r\oplus K_c$ for some Hilbert spaces
$H,K$.

\pf By Corollary 0.7 applied to the identity operator on $E$, if $E$ and $E^*$ are exact,
then $E$ must
be completely isomorphic to a 
``completely complemented" subspace  $F\subset H_r\oplus K_c$,
meaning there is  a $c.b.$  projection $P$ 
from $H_r\oplus K_c$ onto $F$. By [O] this implies that
$E$ is completely isomorphic to  
 $S\oplus T$ for subspaces
$S\subset H_r,  T\subset K_c$. Replacing $H_r,K_c$ by
 these subspaces, we obtain the only if part. The converse is obvious since
$H_r,K_c$ are both exact and (see e.g.  [ER3]) 
$H_r^*=H_c$, $K_c^*=K_r$.\qed

\n {\bf Remark.} The completely isometric analog of the preceding 
characterization is not known:
The only known spaces $E$ such that $ex(E)=ex(E^*)=1$ are
$\ell_\infty^2$, $\ell_1^2$,
$\CC\oplus_\infty R_n$, $\CC\oplus_\infty C_n$ 
as well as their duals,
namely $\CC\oplus_1 C_n$ and $\CC\oplus_1 R_n$ (we refer to [P1] for
a proof that these are exact with constant 1) and
$R_n,C_n$  ($n\ge 0$). Are these the only possible examples ?

% The argument for this is the same as for the preceding corollary.\qed
\medskip

\def\tilde{\widetilde}
Let $E\subset A$ be an operator subspace of a $C^*$-algebra. We will say that
$E$ is {\sl completely complemented\/} if there is a c.b.\ projection $P\colon
\ A\to E$. We have then:

\proclaim Corollary 3.3. Let $A,B$ be $C^*$-algebras. Assume that $E\subset
A$, $E^*\subset B$ (completely isometrically) and that both subspaces are
completely complemented. If either $A$ or $B$ is exact or if $E$ has the CBAP,
then $E$ must be completely isomorphic to $H_r\oplus K_c$ for some Hilbert
spaces $H,K$.

\pf Our assumption implies that the identity on $E$ admits a factorization of
the form
$$E \to A\to E\to B^*\to E$$
hence, by Theorem 0.5, it factors through $H_r\oplus K_c$ and we conclude as
in the preceding corollary using [O].\qed

We now give an application to the operator Hilbert space $OH(I)$
introduced in [P7], to which we refer for more information.

\proclaim Corollary 3.4. Let $A$ be a $C^*$-algebra, let $E\subset A$ be an
exact operator space and let $I$ be an arbitrary set. If a  linear map
$u\colon \ E\to OH(I)$ is c.b.\ then there is a constant $K$ and a state
$f$ on $A$ such that $$\|u(x)\|^2 \le K^2(f(xx^*)f(x^*x))^{1/2}.\leqno
(3.1)\qquad
\forall~x\in E$$
Conversely, for any $E$, any map satisfying (3.1) is  c.b.\ 
More
precisely, we have the estimates
$$K\le 2^{9/4} ex(E)\|u\|_{cb}\quad \hbox{and}\quad \|u\|_{cb} \le
K.$$
Moreover, when $E=A$, the assumption that $E$ is exact is
not needed, and we have the estimates   
$$K\le 2^{9/4}  \|u\|_{cb}\quad \hbox{and}\quad \|u\|_{cb} \le
K.$$

\pf By definition of $OH(I)$,  $u$ is c.b.\ iff the mapping
$\ovl{{}^tu}u\colon \ E\to \ovl{E^*}$ is c.b., and we have (see [P7,
p.~41])
$$\|\ovl{{}^tu}u\|_{cb} = \|u\|^2_{cb}.\leqno
(3.2)$$  Therefore, by
Theorem~0.4, there are states
$f_1,f_2,g_1,g_2$ on $A$ such that $$\|u(x)\|^2 \le 2^{3/2} ex(E)^2
\|u\|^2_{cb} [(f_1(xx^*)g_1(x^*x))^{1/2} + (f_2(x^*x) g_2(xx^*))^{1/2}].\leqno
\forall~x\in E$$ Hence if we let (say) $f = 4^{-1}(f_1+g_1+f_2+g_2)$ we obtain
$$\|u(x)\|^2 \le 2^{9/2} ex(E)^2 \|u\|^2_{cb} (f(xx^*)f(x^*x))^{1/2}.\leqno
\forall~x\in C$$
Thus we obtain the conclusion with $K^2 \le 2^{9/2} ex(E)^2\|u\|^2_{cb}$.

\n Conversely, if (3.1) holds then we have, by Cauchy-Schwarz, for any
$x,y$ in
$E$
$$|\langle u(x),u(y)\rangle| \le \|u(x)\|\  \|u(y)\| \le
K^2(f(xx^*)f(x^*x) f(yy^*)f(y^*y))^{1/4}$$
hence by (1.3)
$$|\langle \ovl{{}^tu}u(x),y\rangle| \le 2^{-1}K^2 (f(xx^*)^{1/2}
f(y^*y)^{1/2} + f(x^*x)^{1/2} f(yy^*)^{1/2}).$$
By the final assertion in Theorem~0.4 and  by (0.8), this implies   
$  \|\ovl{{}^tu}u\|_{cb}\le K^2$,
and by   (3.2) we obtain $\|u\|_{cb} \le K$.
 
\n Now if $E=A$, we may use (ii) in Theorem 0.5 (since that $OH(I)$
has the    CBAP) to justify the last
assertion.
\qed

\n {\bf Remark.} The first part of Corollary 3.4 may fail if 
we do not assume $E$ exact: for instance if $u$ is the identity on
$OH$ and if $T_i$ is an orthonormal basis of $OH$
formed of  self-adjoint operators (see [P7] p. 19), we have
$\sum_1^n \| T_i\|^2 =n$ and $\sum_1^n  (f(T_iT_i^*)f(T_i^*T_i))^{1/2}
=f(\sum_1^n T_i^2)\le \|  \sum_1^n T_i^2  \| = n^{1/2}$.

The next statement improves an unpublished result of Marius Junge (see
[J]) who proved (3.3) (say assuming   $u$ is a complete contraction)
with $(\hbox{\rm Log}(n))^2$ instead of $\hbox{\rm Log}(n)$.
Junge's proof already  used tools from interpolation theory
similar to the ones we use below.

\proclaim Corollary 3.5. Assume $A=B(H)$. Let $u\colon \ A\to OH(I)$ be a
c.b.\ map. Assume (3.1). Then for any $n>1$ and for any $n$-tuple
$x_1,\ldots, x_n$ in A we have
$$\sum^n_1 \|u(x_i)\|^2 \le K^2(c \hbox{ Log}(n) + 1) \left\|\sum^n
x_i\otimes \bar x_i\right\|_{\rm min},\leqno (3.3)$$
where $c>0$ is a numerical absolute constant (independent of $n$ or $u$).

\n {\bf Remark.} If $A$ is an arbitrary $C^*$-algebra, the same argument can
be adapted (using [H2]) to show that (3.3) holds with $\left\|\sum x_i
\otimes
\bar x_i\right\|_{\rm max}$ instead of $\left\|\sum x_i\otimes\bar
x_i\right\|_{\rm min}$.

\n {\bf Note:}\ In the sequel, the constants $c_1,c_2,c_3,\ldots$ will all
be absolute positive numerical constants bounded independently of any
parameter (we can safely say they are all majorized by $10^3$!).

\n To prove Corollary 3.5, we will first need the following consequence of
(3.1).

\proclaim Lemma 3.6. In the situation of Corollary 3.4, assume (for
simplicity) that $A=M_N$. We identify $f$ with an $N\times N$ matrix
$f\ge 0$ with unit trace so that the state $f$ is identified with the map
$x\to \hbox{tr}(fx)$. Then, if (3.1) holds, there are constants $c_1$ and
$c_2$ such that for every $x$ in $E$ and for any $t\ge 2$
$$\|u(x)\|^2 \le K^2\left[c_1 \hbox{\rm Log}(t) \hbox{  tr}(f^{1/2} x
f^{1/2} x^*) + {c_2\over t} f(xx^*+ x^*x)\right]\leqno (3.4)$$

\pf By a change of basis, we may clearly assume that $f$ is a diagonal matrix
so that setting $\lambda_i = f_{ii}$ we have $\lambda_i>0$, $ \sum
\lambda_i=1$ and (3.1) becomes
$$\|u(x)\|^2 \le  K^2\left(\sum_{ij} \lambda_i|x_{ij}|^2 \sum_{ij}
\lambda_j|x_{ij}|^2\right)^{1/2}.\leqno \forall~x\in E$$
Fix a number $t\ge 2$. Let $S(1) = \{(i,j)\mid
t^{-2}\le \lambda_i\lambda^{-1}_j
\le t^2\}$,
$S(2) = \{(i,j)\mid
  \lambda_i\lambda^{-1}_j
> t^2\}$, and 
$S(3) = \{(i,j)\mid
 \lambda_i\lambda^{-1}_j
 <t^{-2} \}$. Then
  let $u_k(x) = u\left(\sum\limits_{ij\in S(k)}
x_{ij}e_{ij}\right)$ ($k=1,2,3$) so that $u
= u_1+u_2+u_3$. Since on $S(2)$ we have $ \lambda_j
\lambda^{-1}_i <t^{-2}$ and similarly on $S(3) $, we can write 
$$\|u_2(x)\|^2 \le  K^2 t^{-1} \sum_{ij} \lambda_i |x_{ij}|^2 = K^2t^{-1}
f(x x^*).\leqno (3.5)$$
and
$$\|u_3(x)\|^2 \le  K^2 t^{-1} \sum_{ij} \lambda_j |x_{ij}|^2 = K^2t^{-1}
f(x^* x).\leqno (3.5)'$$
We now turn to $u_1$. We will use freely the standard notation from
interpolation theory as described e.g.\ in [BL]. Let $E_0$ (resp.\ $E_1$)
be the space $E$ equipped with the norm $f(x^*x)^{1/2}$ (resp.\
$f(xx^*)^{1/2}$). We denote for simplicity $(1\le q<\infty)$ 
$$E(1/2) = (E_0,E_1)_{1/2}\quad \hbox{and}\quad E(1/2,q) =
(E_0,E_1)_{1/2,q}.$$ It is a well known fact that since $E_0,E_1$ are both
Hilbertian we have
$E(1/2) \simeq E(1/2,2)$ with equivalence constants bounded independently of
$N$ or $E$. Moreover, it is a classical fact that
$$\|x\|_{E(1/2)} = \hbox{tr}(f^{1/2}x f^{1/2}x^*).\leqno (3.6)$$
By the extremal property of $E(1/2,1)$ (see [BL, p. 58]),  (3.1) implies
an estimate of the form
$$\|u(x)\| \le c_3K\|x\|_{E(1/2,1)}.~~~~~~~~~~~~~~~~~~~~~~~\leqno
(3.7)\qquad\forall~x\in E$$ Since for any $x$ in $\hbox{span}\{e_{ij}\mid
(i,j)\in S(1)\}$ we have $t^{-1}\|x\|_{E_0} \le \|x\|_{E_1} \le t
\|x\|_{E_0}$, it follows by a well known estimate that $$\|x\|_{E(1/2,1)}
\le c_4(\hbox{Log } t)^{1/2}\|x\|_{E(1/2,2)}$$ hence combining this last
bound with (3.6), (3.7) and the equivalence $E(1/2)  \simeq E(1/2,2)$ we
obtain for any $x$ in
$\hbox{span}\{e_{ij}\mid (i,j)\in S(1)\}$ $$\|u(x)\|^2 \le c_5(\hbox{Log }
t)
\hbox{ tr}(f^{1/2}xf^{1/2}fx^*).$$ A fortiori for all $x$ in $E$ we have
$$\|u_1(x)\|^2 \le c_5(\hbox{Log } t) \hbox{ tr}(f^{1/2}xf^{1/2}x^*).\leqno
(3.8)$$
Now (3.8),  (3.5) and (3.5)' yield (3.4).\qed

\n {\bf Proof of Corollary 3.5.} By (3.4) we have
$$\sum^n_1 \|u(x_i)\|^2 \le K^2\left(c_1 \hbox{ Log}(t)\left\|\sum
x_i\otimes\bar x_i\right\|_{\rm min} + c_2 t^{-1} \left\|\sum
x_ix^*_i\right\|  + c_2t^{-1} \left\|\sum
x^*_ix_i\right\|\right).$$ It is elementary that $\|x_i\|^2 \le
\left\|\sum\limits^n_1 x_i\otimes \bar x_i\right\|$ hence $\left\|\sum
x_ix^*_i\right\|\le \sum \|x_i\|^2 \le n\left\|\sum x_i \otimes\bar
x_i\right\|$, and similarly
$\left\|\sum
x^*_ix_i\right\|\le  n\left\|\sum x_i \otimes\bar
x_i\right\|$. Therefore, choosing $t=n$, we obtain
(3.3).\qed
\medskip

\proclaim Corollary 3.7. Let $E_n\subset B(H)$ be a subspace
$\lambda$-completely isomorphic to $OH_n$, for some $\lambda\ge 1$. Then
for any projection $P\colon\ B(H)\to E_n$ we have
$$\|P\|_{cb}\ge (c_6)^{-1} n^{1/2}(1 + \hbox{\rm Log}(n))^{-1/2} \lambda^{-1}.
\leqno (3.9)$$

\pf Let $v\colon \ E_n\to OH_n$ be an isomorphism such that $\|v\|_{cb}
\|v^{-1}\|_{cb} = \lambda$. Let $u = vP\colon \ B(H)\to OH_n$. Using the
notation of [P6, p.~88], (3.3) implies $\pi^n_{2,oh}(u) \le K(c_1 \hbox{
Log}(n)+1)^{1/2}$ hence by the estimates of $K$ in the last
part of Corollary~3.4
$$\pi^n_{2,oh}(u) \le 2^{9/4}(c_1 \hbox{ Log}(n)+1)^{1/2} \|u\|_{cb} \le
2^{9/4} (c_1 \hbox{ Log}(n) + 1)^{1/2} \|v\|_{cb} \|P\|_{cb}.$$
Let $i\colon \ E\to B(H)$ be the inclusion map. Note $P i=I_E$.
Hence we find 
$$n^{1/2} = \pi^n_{2,oh}(I_{OH_n}) = \pi^n_{2,oh} (uiv^{-1}) \le
\pi^n_{2,oh}(u)\|v^{-1}\|_{cb}\le 2^{9/4} (c_1 \hbox{ Log}(n) + 1)^{1/2}
\lambda\|P\|_{cb}.$$
\qed

\n {\bf Remark.} It remains an open problem (see [P6, Problem 10.2])
whether the logarithmic factor can be entirely removed from either (3.3)
or (3.9).

\medskip

\n {\bf \S 4. Applications to Schur multipliers}

\n Let $B = B(\ell_2)$ and let $K\subset B$ denote the subalgebra of
compact operators. Let $S_1$ denote the trace class, \ie  set of $x$ in
$B$ such that
 ${\rm tr}|x|
< \infty$, where we set$|x|= (x^*x)^{1/2}$, equipped with the trace class
norm
$$\|x\|_1 = {\rm tr}|x|.$$
It is well known that $S^*_1\simeq B$ and $K^*\simeq S_1$ isometrically.
This duality allows us to view  $S_1$ as an operator space for which
the preceding identities become completely isometric (see [ER1] or [P1]).
We will study the Schur multipliers from $K$ to $S_1$ (or from $B$ to
$S_1$), that is to say the linear maps of the form $M_\varphi\colon \
(x_{ij})\to (\varphi_{ij}x_{ij})$ where $(\varphi_{ij})$ is a matrix with
complex entries.

The following is a rather easy consequence of [P4].

\proclaim Theorem 4.1. A Schur multiplier $M_\varphi$ is bounded from $B$ to
$S_1$ (or from $K$ to $S_1$) iff there is a decomposition
$$\varphi_{ij} = a_{ij} + b_{ij}$$
with $\sum\limits_i \sup\limits_j |a_{ij}| + \sum\limits_j \sup\limits_i
|b_{ij}| < \infty$. 

\pf Assume that $M_\varphi\colon \ K\to S_1$ is bounded with norm $\le 1$. For
$a,b$ in $K$, let $u(a,b) = \langle M_\varphi a,b\rangle$. Assuming (say)
$\{a_{ij}\mid i,j\ge 1\}$ finitely supported we have
$$u(a,b) = \sum \varphi_{ij} a_{ij}b_{ij}.$$
By [P4], $u$ satisfies (0.1) and (0.2). As mentioned after (0.2), there
exists a    decomposition   $u = u_1+u_2+u_3+u_4$
so that there are states such that for all $a,b$ in $K$ we have 
$$\eqalign{|u_1(a,b)| &\le K(f_1(aa^*) g_1(b^*b))^{1/2},\quad 
|u_2(a,b)|  \le K(f_2(a^*a) g_2(bb^*))^{1/2}\cr
|u_3(a,b)| &\le K(f_1(aa^*) g_2(bb^*))^{1/2}, 
\quad|u_4(a,b)|  \le K(f_2(a^*a) g_1(b^*b))^{1/2}.}$$
We will use an averaging argument. Let $G$ be the group of all diagonal
unitary matrices on $\ell_2$ equipped with its normalized Haar measure
$m$. 

\n Given a bilinear form $v\colon \ K\times K\to \CC$ we associate to it
the form
$$\tilde v(a,b) = \int\limits_{G\times G} v(w_1aw_2, w^{-1}_2bw^{-1}_1)~~
dm(w_1)dm(w_2).$$
Clearly $\|\tilde v\| \le \|v\|$. Moreover, $\tilde v(w_1aw_2,
w^{-1}_1bw^{-1}_2) = \tilde v(a,b)$ for any $w_1,w_2$ in $G$. Therefore there
is a Schur multiplier $\psi$ such that, if $a,b$ are finitely supported, we
have
$$\tilde v(a,b) = \sum \psi_{ij}a_{ij}b_{ij}.$$
We now apply this averaging procedure to each of $u_1,u_2,u_3$ and $u_4$. For
$u_1$ we obtain (by Cauchy-Schwarz):
$$\eqalign{|\tilde u_1(a,b)| &\le \int|u(w_1aw_2,w^{-1}_1bw^{-1}_2)|~~ dm(w_1)
dm(w_2)\cr
&\le K\left(\int f_1(w_1aa^*w^{-1}_1) dm(w_1)\right)^{1/2} \left(\int
g_1(w_2b^*bw^{-1}_2)dm(w_2)\right)^{1/2}.}$$
Moreover $\int w^{-1}_1f_1w_1\ dm(w_1)$ and $\int w^{-1}_2g_1w_2\ dm(w_2)$ are
diagonal states on $K$.  Setting $x^1_i = (f_1)_{ii}$ and $y^1_j = (g_1)_{jj}$
and $\varphi^1_{ij} = \tilde u_1(e_{ij},e_{ij})$ we find
$$\left|\sum\varphi^1_{ij} a_{ij}b_{ij}\right| \le K\left(\sum_{ij}
x^1_i|a_{ij}|^2 \cdot \sum_{ij} |a_{ij}|^2 y^1_j\right)^{1/2}.$$
Hence $|\varphi^1_{ij}| \le K|x^1_i|^{1/2} |y^1_j|^{1/2}$ with
$$\sum |x^1_i|= 1,\quad \sum |y^1_j| = 1.$$
Applying the same procedure successively to $u_2,u_3$ and $u_4$, we obtain a
decomposition $\varphi = \varphi^1 + \varphi^2 + \varphi^3 + \varphi^4$ where
$\varphi^2$ satisfies (by symmetry) a similar bound as $\varphi^1$ and
$\varphi^3,\varphi^4$ are such that there are $x^3_i,y^3_i$, $x^4_i,y^4_i$ all
$\ge 0$ and with $\sum x^3_i = \sum y^3_i = \sum x^4_i = \sum y^4_i = 1$ such
that $|\varphi^3_{ij}|\le K(x^3_iy^3_i)^{1/2}$ and $|\varphi^4_{ij}| \le
K(x^4_jy^4_j)^{1/2}$. Clearly we have $\sum\limits_i \sup\limits_j
|\varphi^3_{ij}|\le K$ and $\sum\limits_j \sup\limits_i |\varphi^4_{ij}|\le K$.

\n On the other hand, we claim that $\varphi^1$ and $\varphi^2$ each admit
a decomposition of the kind described in Theorem~4.1. Indeed, it suffices
to prove this with a uniform bound when $\varphi^1,x^1,y^1$ are all
finitely supported. Then, after a suitable permutation of the indices, we
may assume that $x^1_i$ and $y^1_j$ are both non-decreasing. Let then
$\alpha_{ij} =
\varphi^1_{ij}$ if $i\le j$ and $\beta_{ij} = \varphi^1_{ij}$ if $i>j$. We
have then $\alpha_{ij} \le (x^1_iy^1_j)^{1/2}\le (x^1_jy^1_j)^{1/2}$ and
$\beta_{ij} \le (x^1_iy^1_j)^{1/2}\le (x^1_iy^1_i)^{1/2}$ hence we obtain
$\varphi^1_{ij} = \alpha_{ij} +\beta_{ij}$ with
$$\sum_j \sup_i |\alpha_{ij}| \le 1\quad \hbox{and}\quad \sum_i \sup_j
|\beta_{ij}|\le 1,$$
which proves our claim.
Using this we finally obtain a decomposition $\varphi_{ij} = a_{ij} + b_{ij}$
with $\sum\limits_i \sup\limits_j |a_{ij}|\le K+2$ and $\sum\limits_j
\sup\limits_i |b_{ij}| \le K+2$. This proves the ``only if'' part.

\n  Conversely,
if $\varphi$ can be decomposed as in Theorem~4.1, we will show that
$M_\varphi\colon \ B\to S_1$ is bounded. It clearly suffices to show that
$\sum\limits_i \sup\limits_j |\varphi_{ij}|<\infty$ and $\sum\limits_j
\sup\limits_i |\varphi_{ij}|<\infty$ are sufficient conditions for the
boundedness of $M_\varphi\colon \ B\to S_1$. But this is obvious since:
$$\left|\sum \varphi_{ij}a_{ij}b_{ij}\right|  \le \sum_i \sup_j
|\varphi_{ij}| \sum_j |a_{ij}b_{ij}|
 \le \sum_i \sup_j |\varphi_{ij}|\ \|a\|_B\ \|b\|_B$$
 {and similarly}
$$\left|\sum \varphi_{ij}a_{ij}b_{ij}\right| \le \sum_j \sup_i
|\varphi_{ij}|\
\|a\|_B\ \|b\|_B.$$
\qed

We now turn to the c.b.\ analogue of the preceding statement.

\proclaim Theorem 4.2. Consider a complex matrix $(\varphi_{ij})$. The
following assertions are equivalent:
\item{(i)} $M_\varphi$ is c.b.\ from $K$ to $S_1$.
\item{(i)$'$} $M_\varphi$ is c.b.\ from $B$ to $S_1$.
\item{(ii)} There are $x,y$ in $\ell_2$ and a constant $C$ such that
$|\varphi_{ij}|\le C|x_i|\ |y_j|$  for all~$i,j$.
\item{(iii)} There is a element $T$ in $S_1$ and a constant $C$ such that
$|\varphi_{ij}| \le C|T_{ij}|$ forall~$i,j$.

\pf The main implication is (i) $\Rightarrow$ (ii). The fact that (ii)
$\Leftrightarrow$ (iii) is elementary; we include it for the record. Assume
(i). Let $u\colon \ K\times K\to \CC$ be as before  with $\|u\|_{cb}\le
1$. By Theorem~0.5, there is a decomposition $u = u_1+u_2$ and states
$f_1,f_2$, $g_1,g_2$ such that, for some constant $C$ we have for all
$a,b$ in $K$ :
$$\eqalign{|u_1(a,b)| &\le C(f_1(aa^*) g_1(b^*b))^{1/2}\cr
|u_2(a,b)| &\le C(f_2(a^*a) g_2(bb^*))^{1/2}.}$$
Applying the same averaging procedure as above, we find a decomposition $u=
\tilde  u = \tilde u_1+\tilde u_2$ and matrices $\varphi^1$ and $\varphi^2$
such that $\tilde u_1 = M_{\varphi^1}$, and $\tilde u_2 = M_{\varphi^2}$.
Moreover, we have
$$\eqalign{\left|\sum \varphi^1_{ij}a_{ij}b_{ij}\right| &\le C\left(\sum_{ij}
x^1_i|a_{ij}|^2 \sum_{ij} y^1_j|a_{ij}|^2\right)^{1/2}\cr
\left|\sum \varphi^2_{ij}a_{ij}b_{ij}\right| &\le C\left(\sum_{ij} x^2_j
|a_{ij}|^2 \sum_{ij} y^2_i |a_{ij}|^2\right)^{1/2}}$$
where $x^1_i = (f_1)_{ii}$, $x^2_j = (f_2)_{jj}$, $y^1_j = (g_1)_{jj}$ and
$y^2_i = (g_2)_{ii}$. This gives us
$$|\varphi_{ij}| = |\varphi^1_{ij}+\varphi^2_{ij}|\le C(x^1_iy^1_j +
x^2_jy^2_i)^{1/2}$$ hence if we set $x_i = (|x^1_i| + |y^2_i|)^{1/2}$ and $y_j
= (|y^1_j| + |x^2_j|)^{1/2}$, we obtain
$$|\varphi_{ij}| \le Cx_iy_j.$$
This shows that (i) $\Rightarrow$ (ii). (ii) $\Rightarrow$ (iii) is obvious.
The converse is elementary. Indeed, if $T = \sum\limits_k x^k\otimes y^k
\lambda_k$ with $x^k,y^k$ in the unit ball of $\ell_2$ and $\sum|\lambda_k| <
\infty$, we have
$$|T_{ij}|\le X_iY_j\quad \hbox{where}\quad X_i = \left(\sum_k \lambda_k
|x^k_i|^2\right)^{1/2}$$
and $Y_j = \left(\sum\limits_k \lambda_k |y^k_j|^2\right)^{1/2}$.

\n It remains only to show (ii) $\Rightarrow$ (i)$'$ (since (i)$'
\Rightarrow$ (i) is trivial by restriction). Assume (ii). Let $f = \sum
|x_i|^2e_{ii}$ and
$g = \sum|y_j|^2 e_{jj}$. For any finitely supported $a,b$ we have by
Cauchy-Schwarz
$$\eqalign{\left|\sum \varphi_{ij}a_{ij}b_{ij}\right| &\le C\sum |x_i|\
|a_{ij}|\ |y_j| \le C\left(\sum_{ij}|x_i|^2 |a_{ij}|^2 \sum_{ij} |a_{ij}|^2
|y_j|^2\right)^{1/2}\cr
&\le C(f(aa^*) g(b^*b))^{1/2}.}$$
Hence we conclude that the bilinear from $u$ is c.b.; {\sl a fortiori\/} $u$
is j.c.b.\ and $M_\varphi$ is c.b.\qed 
\medskip

\n {\bf Remark.} The preceding argument actually shows that the bilinear form
$(a,b)\to \langle M_\varphi a,b\rangle$ is j.c.b.\ iff it is c.b.\medskip

\n {\bf Remark.} It is easy to see that (ii) in Theorem~4.2 is equivalent to:
\item{(ii)$'$} There are $(a_{ij})$ and $(b_{ij})$ with $\sum\limits_i
\sup\limits_j |a_{ij}| + \sum\limits_j \sup\limits_i |b_{ij}| < \infty$ such
that
$$|\varphi_{ij}| = \sqrt{a_{ij}b_{ij}}.$$

Thus, whereas we had the arithmetic mean in Theorem~4.1, here we find the
geometric mean of the same terms.

\n In particular, it is easy to see from this that there are bounded Schur
multipliers $M_\varphi\colon \ K\to S_1$ which are {\sl not\/} c.b.

 \bigskip

\centerline {\bf References}

 \item{[B1]} D. Blecher. Generalizing
Grothendieck's program.  "Function spaces",
Edited by K. Jarosz, Lecture Notes in Pure and
Applied Math. vol.136, Marcel Dekker, 1992.

\item{[BL]} J. Bergh and J. L\"ofstr\"om. Interpolation spaces. An
 introduction. Springer Verlag, New York. 1976.

 \item{[BP]} D. Blecher and V. Paulsen. Tensor products of operator spaces.
 J. Funct. Anal. 99 (1991) 262-292.
 
 \item{[C]} A.~Connes,  {Almost periodic states and factors of type {III$_1$}}, J.
  Funct. Anal.  {16} (1974), 415--455.

\item{[CES]}  E.  Christensen, E. Effros and A. Sinclair. 
Completely bounded multilinear maps and
$C^*$-algebraic cohomology.
Invent. Math. 90 (1987) 279-296.

\item{[CS1]} E. Christensen and A. Sinclair. 
Representations of completely bounded multilinear operators.
J. Funct. Anal. 72 (1987) 151-181.

\item{[CS2]} $\underline{\hskip1.5in}$.   A survey of
completely bounded operators.	Bull. London Math. Soc.
21 (1989) 417-448.

 \item{[D]} 
K.~Dykema,  {Free products of finite dimensional and other von {Neumann}
  algebras with respect to non-tracial states}, Free Probability (D.-V.
  Voiculescu, ed.), Fields Institute Communications, vol.~12, American
  Mathematical Society, 1997, pp.~41--88.

\item{[ER1]} E. Effros and Z.J. Ruan. Operator
Spaces.  Oxford Univ. Press,
Oxford, 2000.

\item{[ER2]} $\underline{\hskip1.5in}$.  A new
approach to operator spaces.
 Canadian Math. Bull.
34 (1991) 329-337.

\item{[ER3]} $\underline{\hskip1.5in}$. Self duality for the Haagerup
 tensor product and Hilbert space factorization.  J. Funct. Anal. 100
(1991) 257-284.

\item{[G]} A. Grothendieck.  R\'esum\'e de la th\'eorie
m\'etrique des produits tensoriels topologiques.  Boll..
Soc. Mat. S$\tilde{a}$o-Paulo  8 (1956), 1-79.

   \item{[H1]} U. Haagerup. 
The Grothendieck inequality for bilinear forms on $C\sp *$-algebras. Adv. in
Math. 56 (1985)  93--116.

\item{[H2]} $\underline{\hskip1.5in}$. Self-polar forms,
 conditional expectations and the weak expectation property
for $C^*$-algebras. Unpublished manuscript (1995).

\item{[HI]} U. Haagerup and  T. Itoh. Grothendieck 
type norms for bilinear forms on $C\sp *$-algebras.
J. Operator Theory 34 (1995)  263--283.
 
\item{[J]} M. Junge. A first attempt to the little Grothendieck inequality and
related embeddings. Unpublished notes from a  lecture at I.H.P.\ (January 2000).

\item{[JP]} M. Junge and G. Pisier. Bilinear
 forms on exact operator
spaces and $B(H)\otimes B(H)$. Geometric and 
Functional Analysis (GAFA Journal) 
5 (1995) 329-363.
 
 \item{[Ki1]} E. Kirchberg.  On non-semisplit
extensions, tensor products and exactness of
group $C^*$-algebras. Invent. Math. 112 (1993)
449-489.

\item{[Ki2]} $\underline{\hskip1.5in}$.
 Commutants of unitaries in UHF algebras and functorial properties of
exactness. J. Reine Angew. Math. 452 (1994), 39--77.

\item{[Ki3]} $\underline{\hskip1.5in}$.
 Exact $C^*$-algebras, Tensor products, and Classification
of purely infinite algebras. Proceedings ICM 94, Z\"urich,
vol. 2, 943-954, Birkhauser, 1995.

\item{[L]} E.C. Lance.  On nuclear $C\sp{*} $-algebras.
 J. Functional Analysis 12 (1973), 157--176. 

%\item{[LPP]} F. Lust-Piquard and G. Pisier. Non-commutative Khintchine and Paley
%inequalities. Ark. f\"or Mat. 29 (1991) 241-260.

\item{[O]} T. Oikhberg. Direct sums of homogeneous operator spaces.
J. London Math. Soc. To appear.

 \item{[OP]} T. Oikhberg and G. Pisier.
 The ``maximal" tensor product of
two operator spaces. Proc. Edinburgh Math. Soc. 42 (1999)
267--284.

 \item{[PaS]} V. Paulsen and R. Smith. Multilinear maps
and tensor norms on operator systems. J. Funct. Anal. 73
(1987) 258-276.

\item{[P1]} G. Pisier.  Introduction to     operator space theory.
Cambridge Univ. Press. Book to appear.

\item{[P2]} $\underline{\hskip1.5in}$. A simple proof
of a theorem of Kirchberg
and related results on
$C^*$-norms.
J. Op. Theory.  35 (1996) 317-335.

\item{[P3]} $\underline{\hskip1.5in}$. Exact operator
spaces. Colloque sur les alg\`ebres d'op\'era\-teurs.
in ``Recent advances in operator algebras" (Orl\'eans 1992)
 Ast\'erisque (Soc. Math. France) 232 (1995) 159-186.

\item{[P4]} $\underline{\hskip1.5in}$. Grothendieck's theorem for
non-commutative C$^*$-algebras with an appendix on \-Grothendieck's
constants. J. Funct. Anal. 29 (1978) 397-415.

\item {[P5]} $\underline{\hskip1.5in}$. Factorization of linear
operators and the Geometry of Banach spaces.  CBMS
(Regional conferences of the A.M.S.)    60, (1986),
Reprinted with corrections 1987.

\item {[P6]} $\underline{\hskip1.5in}$.  Non-commutative
 vector valued $L_p$-spaces and completely
$p$-summing maps. Soc. Math. France. 
Ast\'erisque 247 (1998)  1-131. 

\item {[P7]} $\underline{\hskip1.5in}$. 
The operator Hilbert space $OH$,
complex interpolation and tensor norms.
Memoirs Amer.
Math. Soc.  vol. 122 , 585 (1996) 1-103.

\item{[Sc]} R. Schneider.  Zonoids whose polars are zonoids. Proc.
Amer. Math. Soc. 50 (1975) 365--368. 

\item{[S1]}  D. Shlyakhtenko. Free
quasi-free states. Pacific J. Math. 177 (1997),
no. 2, 329--368.

 \item{[S2]} $\underline{\hskip1.5in}$.  Limit distributions
of matrices with bosonic and fermionic entries.
Free probability theory (Waterloo, ON, 1995),
241--252, Fields Inst. Commun., 12, Amer. Math.
Soc., Providence, RI, 1997. 

 \item{[S3]} $\underline{\hskip1.5in}$.   $A$-valued semicircular systems. J. Funct. Anal.
166 (1999),   1--47.

\item{[VDN]} D. Voiculescu, K. Dykema, A. Nica. Free
random variables. CRM Monograph Series,
Vol. 1,  Amer. Math. Soc., Providence RI.

\end